\documentclass[10pt]{article}

\usepackage{graphicx}

\usepackage{amssymb} 
\usepackage{amsmath}
\usepackage{times}
\usepackage{bm}
\newcommand{\vvec}{\textrm{vec}}
\newcommand{\tr}{\textrm{tr}}

\newcommand{\bfsym}{\boldsymbol}  

\title{Improved estimation in a general multivariate elliptical model}
\author
{
   Tatiane F. N. Melo\\
     {\it \footnotesize Institute of Mathematics and Statistics, Federal University of Goi\'as, Brazil} \vspace{-0.2cm}\\
     {\footnotesize email: {\tt tmelo@ufg.br}} \\ \\
   Silvia L. P. Ferrari\\
   {\it \footnotesize Department of Statistics, University of S\~ao Paulo, Brazil} \vspace{-0.2cm}\\
     {\footnotesize email: {\tt silviaferrari@usp.br}} \\ \\
   Alexandre G. Patriota\\
     {\it \footnotesize Department of Statistics, University of S\~ao Paulo, Brazil} \vspace{-0.2cm}\\
     {\footnotesize email: {\tt patriota@ime.usp.br}} \\ \\
}
\date{}

\begin{document}

\maketitle

\begin{abstract}
The problem of reducing the bias of maximum likelihood estimator in a general multivariate elliptical regression model is considered. The model is very flexible and allows the mean vector and the dispersion matrix to have parameters in common. Many frequently used models are special cases of this general formulation, namely: errors-in-variables models, nonlinear mixed-effects models, heteroscedastic nonlinear models, among others. In any of these models, the vector of the errors may have any multivariate elliptical distribution. We obtain the second-order bias of the maximum likelihood estimator, a bias-corrected estimator, and a bias-reduced estimator. Simulation results indicate the effectiveness of the bias correction and bias reduction schemes.

\end{abstract}


\section{Introduction}\label{introduction}  

It is well known that, under some standard regularity conditions, maximum-like\-li\-hood estimators (MLEs) are consistent and asymptotically normally distributed. Hence, their biases converge to zero when the sample size increases. However, for finite sample sizes, the MLEs are in general biased and bias correction plays an important role in the point estimation theory.

A general expression for the term of order $O(n^{-1})$ in the expansion of the bias of MLEs was given by Cox and Snell (1968). This term is often called second-order bias and can be useful in actual problems. For instance, a very high second-order bias indicates that other than maximum-likelihood estimation procedures should be used. 
Also, corrected estimators can be formulated by subtracting the estimated second-order biases from the respective MLEs. 
It is expected that these corrected estimators have smaller biases than the uncorrected ones, especially in small samples.
 
Cox and Snell's formulae for second-order biases of MLEs were applied in many models. Cordeiro and McCullagh (1991) use these formulae in generalized linear models; Cordeiro and Klein (1994) compute them for ARMA models; Cordeiro et al. (2000) apply them for symmetric nonlinear regression models; Vasconcellos and Cordeiro (2000) obtain them for multivariate nonlinear Student t regression models. More recently, Cysneiros et al. (2010) study the univariate heteroscedastic symmetric nonlinear regression models (which are an extension of Cordeiro et al. 2000) and Patriota and Lemonte (2009) obtain a general matrix formula for the bias correction in a multivariate normal model where the mean and the covariance matrix have parameters in common.

An alternative approach to bias correction was suggested by Firth (1993). The idea is to adjust the estimating function so that the estimate becomes less biased. This approach can be viewed as a ``preventive" method, since it modifies the original score function, prior to obtaining the parameter estimates. In this paper, estimates obtained from Cox and Snell's approach and Firth's method will be called bias-corrected estimates and bias-reduced estimates, respectively. Firth showed that in generalized linear models with canonical link function the preventive method is equivalent to maximizing a penalized likelihood that is easily implemented via an iterative adjustment of the data. The bias reduction proposed by Firth has received considerable attention in the statistical literature. For models for binary data, see Mehrabi and Matthews (1995); for censored data with exponential lifetimes, see Pettitt et al. (1998). In Bull et al. (2002) bias reduction is obtained for the multinomial logistic regression model. In Kosmidis and Firth (2009) a family of bias-reducing adjustments was developed for a general class of univariate and multivariate generalized nonlinear models. The bias reduction in cumulative link models for ordinal data was studied in Kosmidis (2014). Additionally, Kosmidis and Firth (2011) showed how to obtain
the bias-reducing penalized maximum likelihood estimator by using the equivalent Poisson log-linear model for the parameters of a multinomial logistic regression.

It is well-known and was noted by Firth (1993) and Kosmidis and Firth (2009) that the reduction in bias may sometimes be accompanied by inflation of variance, possibly yielding an estimator whose mean squared error is bigger than that of the original one. Nevertheless, published empirical studies such as those mentioned above show that, in some frequently used models, bias-reduced and bias-corrected estimators can perform better than the unadjusted maximum likelihood estimators, especially when the sample size is small.

Our goal in this paper is to obtain bias correction and bias reduction to the maximum likelihood estimators for the general multivariate elliptical model. We extend the work of Patriota and Lemonte (2009) to the elliptical class of distributions defined in Lemonte and Patriota (2011). We focus on analytical methods only, because simulations for a general multivariate normal model suggests that analytical bias corrections outperforms the computationally intensive bootstrap methods (Lemonte, 2011).

In order to illustrate the ampleness of the general multivariate elliptical model, we mention some of its submodels: multiple linear regression, heteroscedastic multivariate nonlinear regressions, nonlinear mixed-effects models (Patriota 2011), heteroscedastic errors-in-variables models (Patriota et al. 2009a,b), structural equation models, multivariate normal regression model with general parametrization (Lemonte, 2011), simultaneous equation models and mixtures of them. It is important to note that the usual normality assumption of the error is relaxed and replaced by the assumption of elliptical errors. The elliptical family of distributions includes many important distributions such as multivariate normal, Student $t$, power exponential, contaminated normal, Pearson II, Pearson VII, and logistic, with heavier or lighter tails than the normal distribution; see Fang et al. (1990).

The paper is organized as follows. Section \ref{BiasCorrectionReduction} presents the notation and general results for bias correction and bias reduction. Section \ref{model} presents the model and our main results, namely the general expression for the second-order bias of MLEs, in the general multivariate elliptical model. Section \ref{special-cases} applies our results in four important special cases: heteroscedastic nonlinear (linear) model, nonlinear mixed-effects models, multivariate errors-in-variables models and log-symmetric regression models. Simulations are presented in Section \ref{simulation}. Applications that use real data are presented and discussed in Section \ref{applications}. Finally, Section \ref{conclusion} concludes the paper. Technical details are collected in one appendix.


\section{Bias correction and bias reduction}\label{BiasCorrectionReduction}

Let $\theta$ be the $p$-vector of unknown parameters and ${\theta}_{r}$ its $r$th element. Also, let $U(\theta)$ be the score function and $U_r(\theta) = U_r$ its $r$th element. We use the following tensor notation for the cumulants of the log-likelihood derivatives introduced by Lawley
(1956):

\[\kappa_{rs} = E\bigg(\frac{\partial U_r}{\partial\theta_{s}}\bigg),\quad \kappa_{r,s} = E(U_r U_s),\quad 
\kappa_{rs,t} = E\bigg(\frac{\partial U_r}{\partial\theta_{s}} U_{t}\bigg),\]
\[\kappa_{rst} = E\bigg(\frac{\partial^{2} U_r}{\partial\theta_{s}\partial\theta_{t}}\bigg),\quad
\kappa_{rs}^{(t)} = \frac{\partial\kappa_{rs}}{\partial\theta_{t}}, \quad \kappa_{r,s,t} = E(U_r U_s U_t),\] 
and so on. The indices $r$, $s$ and $t$ vary from $1$ to $p$. The typical $(r,s)$th element of the Fisher information matrix $K(\theta)$ is $\kappa_{r,s}$ and we denote by $\kappa^{r,s}$ the corresponding element of $K(\theta)^{-1}$. All $\kappa$'s refer to a total over the sample and are, in general, of order $n$. Under standard regular conditions, we have that $\kappa_{rs} = - \kappa_{r,s}$, $\kappa_{rs,t} = \kappa_{rs}^{(t)}- \kappa_{rst} $ and $\kappa_{r,s,t} = 2\kappa_{rst} -\kappa_{rs}^{(t)} - \kappa_{rt}^{(s)} - \kappa_{st}^{(r)}$. These identities will be used to facilitate some algebraic operations.

Let $B_{\widehat{\theta}}(\theta)$ be the second-order bias vector of $\widehat{\theta}$ whose $j$th element is $B_{\widehat{\theta}_{j}}(\theta)$, $j=1,2,\ldots,$ $p$. It follows from the general expression for the multiparameter second-order biases of MLEs given by Cox and Snell (1968) that
\begin{equation}\label{biascorrection}
B_{\widehat{\theta}_{j}}(\theta) = \sum_{r,s,t = 1}^p\kappa^{j,r}\kappa^{s,t}\biggl\{
\frac{1}{2}\kappa_{rst}+\kappa_{rs,t}\biggr\}. 
\end{equation}

The bias corrected MLE is defined as 
\begin{equation*}{\widehat{\theta}}_{BC} = \widehat{{\theta}} - B_{\widehat{\theta}}(\widehat{\theta}).\end{equation*}
The bias-corrected estimator ${\widehat{\theta}}_{BC}$ is expected to have smaller bias than the uncorrected estimator, $\widehat{{\theta}}$.

Firth (1993) proposed an alternative method to partially remove the bias of MLEs. The method replaces the score function by its modified version
\begin{equation*}\label{scoreModified}
U^*(\theta) = U(\theta) - K(\theta) B_{\widehat{\theta}}(\theta),
\end{equation*}
and a modified estimate, ${\widehat{\theta}}_{BR}$, is given as a solution to $U^*(\theta) = 0$. It is noticeable that, unlike Cox and Snell's approach, Firth's bias reduction method does not depend on the finiteness of ${\widehat{\theta}}$. 

\section{Model and main results}\label{model}

We shall follow the same notation presented in Lemonte and Patriota (2011). The elliptical model as defined in Fang et al. (1990) follows. A $q\times 1$ random vector ${Y}$ has a multivariate elliptical distribution with location parameter ${\mu}$ and a definite positive scale matrix
${\Sigma}$ if its density function is 
\begin{equation}\label{dens}
f_{{Y}}({y}) = |{\Sigma}|^{-1/2} g\bigl(({y} - {\mu})^{\top}{\Sigma}^{-1}({y} - {\mu})\bigr),
\end{equation}
where $g:[0,\infty)\to(0,\infty)$ is called the density generating function, and it is such that $\int_{0}^{\infty}u^{\frac{q}{2}-1}g(u) du < \infty$. We will denote ${Y}\sim El_{q}({\mu},{\Sigma}, g) \equiv El_{q}({\mu},{\Sigma})$. It is possible to show that the characteristic function is $\psi(t) = \mbox{E}(\exp(it^\top Y))=\exp(it^\top\mu)\varphi(t^\top\Sigma t)$, where $t \in \mathbb{R}^q$ and $\varphi: [0,\infty) \to \mathbb{R}$. Then, if $\varphi$ is twice differentiable at zero, we have that $\mbox{E}(Y) = \mu$ and $\mbox{Var}(Y) = \xi \Sigma$, where $\xi = \varphi'(0)$. We assume that the density generating function $g$ does not have any unknown parameter, which implies that $\xi$ is a known constant. From (\ref{dens}), when ${\mu} = {0}$ and ${\Sigma} = {I}_{q}$, where ${I}_{q}$ is a $q\times q$ identity matrix, we obtain the spherical family of densities. A comprehensive exposition of the elliptical multivariate class  of distributions can be found in Fang et al. (1990). Table 1 presents the density generating functions of some multivariate elliptical distributions.

\begin{table*}
\caption{Generating functions of some multivariate elliptical distributions.}
\label{Tab1}
\begin{center}
\begin{tabular}{cc}\hline \noalign{\smallskip}\\
   Distribution     & Generating function $g(u)$  \\
   \noalign{\smallskip}\hline\noalign{\smallskip}
      normal        & \!\!\!\!\!\!\!\!\!\!\!\!\!\!\!\!\!\!\!\!\!\!\!\!\!\!\!\!\!\!\!\!\!\!\!\!\!\!\!\!\!\!\!\!\!\!\!\!\!\!\!\!\!\!\! $ \frac{1}{(\sqrt{2 \pi})^q} \ e^{-u/2}$ \\
 \\
      Cauchy        & \!\!\!\!\!\!\!\!\!\!\!\!$ \frac{\Gamma\left(\frac{1+q}{2}\right)}{\Gamma\left(\frac{1}{2}\right)} \pi^{-q/2} (1 + u)^{-(1 + q)/2}$  \\
 \\
     Student $t$      & \ \ \ \ \ \ \ \ \ \ \ \ \ \ \ \ $ \frac{\Gamma\left(\frac{\nu + q}{2}\right)}{\Gamma\left(\frac{\nu}{2}\right)} \pi^{-q/2} \nu^{-q/2} \left(1 + \frac{u}{\nu}\right)^{-(\nu + q)/2}$, \ $\nu > 0$ \\
\\
power exponential   & $ \frac{\lambda \Gamma\left(\frac{q}{2}\right)}{\Gamma\left(\frac{q}{2\lambda}\right)} 2^{-q/(2 \lambda)} \pi^{-q/2} e^{-u^{\lambda}/2}$, \ $\lambda > 0$\\
\noalign{\smallskip}\hline 
\end{tabular}
\end{center}
\end{table*}

Let $Y_1 , Y_2 , . . . , Y_n$ be $n$ independent random vectors, where $Y_i$ has dimension $q_i\in \mathbb{N}$, for $i=1,2, . . ., n$. The general multivariate elliptical model (Lemonte and Patriota 2011) assumes that 
\begin{equation*} 
{Y}_{i} = {\mu}_{i}({\theta},{x}_{i}) + {e}_{i},\quad i=1,\ldots,n,
\end{equation*}
with ${e}_{i} \stackrel{ind}{\sim} El_{q_{i}}({0},{\Sigma}_i({\theta}, {w}_{i}))$, where ``$ \stackrel{ind}{\sim}$'' means ``independently distributed as'', ${x}_i$ and ${w}_i$ are $m_i\times 1$ and $k_{i}\times 1$ nonstochastic vectors of auxiliary variables, respectively, associated with the $i$th observed response ${Y}_{i}$, which may have components in common. Then,

\begin{equation}\label{MainModel}
{Y}_{i} \stackrel{ind}{\sim}El_{q_{i}}({\mu}_{i}, {\Sigma}_i),\quad i=1,\ldots,n,
\end{equation}

\noindent where ${\mu}_{i} = {\mu}_{i}({\theta},{x}_{i})$ is the location parameter and ${\Sigma}_i = {\Sigma}_i({\theta}, {w}_{i})$ is the definite positive scale matrix. Both ${\mu}_{i}$ and ${\Sigma}_i$ have known functional forms and are twice differentiable with respect to each element of ${\theta}$. Additionally, ${\theta}$ is a $p$-vector of unknown parameters (where $p<n$ and it is fixed). Since ${\theta}$ must be identifiable in model (\ref{MainModel}), the functions ${\mu}_i$ and ${\Sigma}_i$ must be defined to accomplish such restriction. 

Several important statistical models are special cases of the general formulation (\ref{MainModel}), for example, linear and nonlinear regression models, homoscedastic or heteroscedastic measurement error models, and mixed-effects models with normal errors. It is noteworthy that the normality assumption for the errors may be relaxed and replaced by any distribution within the class of elliptical distributions, such as the Student $t$ and the power exponential distributions. The general formulation allows a wide range of different specifications for the location and the scale parameters, coupled with a large collection of distributions for the errors. Section 4 presents four important particular cases of the main model (\ref{MainModel}) that show the applicability of the general formulation.

For the sake of simplifying the notation, let ${z}_i = {Y}_i - {\mu}_i$ and $u_{i} = {z}_i^{\top}{\Sigma}_i^{-1}{z}_i$. The log-likelihood function associated with (\ref{MainModel}), is given by 
\begin{equation} \label{log-likelihood}
\ell({\theta}) = \sum_{i=1}^n \ell_{i}({\theta}),
\end{equation}
where $\ell_{i}({\theta}) = -\frac{1}{2} \log{|{\Sigma}_i|} + \log g(u_i)$. It is assumed that $g(\cdot)$, ${\mu}_i$ and ${\Sigma}_i$ are such that $\ell({\theta})$ is a regular log-likelihood function (Cox and Hinkley 1974, Ch. 9) with respect to ${\theta}$. To obtain the score function and the Fisher information matrix, we need to derive $\ell(\theta)$ with respect to the unknown parameters and to compute some moments of such derivatives. We assume that such derivatives exist. Thus, we define

\[
{a}_{i(r)} = \frac{\partial {\mu}_i}{\partial \theta_{r}}, \quad
{a}_{i(sr)} = \frac{\partial^2 {\mu}_i}{\partial
  \theta_{s} \partial\theta_{r}}, \quad {C}_{i(r)} = \frac{\partial
  {\Sigma}_i}{\partial \theta_{r}}, \quad {C}_{i(sr)} =
\frac{\partial^2 {\Sigma}_i}{\partial \theta_{s} \partial\theta_{r}}
\]
and\[
{A}_{i(r)} = -{\Sigma}_i^{-1} {C}_{i(r)}{\Sigma}_i^{-1},
\]
for $r,s = 1,\ldots,p$. We make use of matrix differentiation methods (Magnus and Neudecker 2007) to compute the derivatives of the log-likelihood function. The score vector and the Fisher information matrix for ${\theta}$ can be shortly written as 
\begin{equation}\label{ScoreFisher2}
U(\theta) = {F}^{\top}{H}{s} \quad \mbox{and}\quad  K(\theta) ={F}^{\top}\widetilde{H}{F},
\end{equation}
respectively, with $F = \left(F_1^\top, \ldots, F_n^\top\right)^{\top}$, $H = \mbox{block-diag}\left\{H_1, \ldots, H_n\right\}$, $s = ( s_1^\top, \ldots, s_n^\top )^{\top}$, $\widetilde{H} = H M H$ and $M = \mbox{block-diag}\left\{ M_1^\top,\ldots, M_n^\top \right\}$, wherein
\begin{equation*} \label{matrix-score}
{F}_i =
  \begin{pmatrix}
    {D}_i\\
    {V}_i\\
  \end{pmatrix},\quad
{H}_i =
  \begin{bmatrix}
    {\Sigma}_i & {0}\\
    {0} & 2{\Sigma}_i\otimes {\Sigma}_i
  \end{bmatrix}^{-1}, \quad
{s}_i =
  \begin{bmatrix}
    v_{i}{z}_i\\
    -\vvec({\Sigma}_i - v_{i}{z}_i{z}_i^{\top})
  \end{bmatrix},
\end{equation*}
where the ``vec" operator transforms a matrix into a vector by stacking the columns of the matrix, ${D}_i = ({a}_{i(1)}, \ldots, {a}_{i(p)})$, ${V}_i = (\vvec({C}_{i(1)}),\ldots, \vvec({C}_{i(p)}))$, $v_{i} = -2W_{g}(u_{i})$ and $W_g(u) = \mbox{d} \log g(u)/\mbox{d} u$. 
Here, we assume that $F$ has rank $p$ (i.e., ${\mu}_i$ and ${\Sigma}_i$ must be defined to hold such condition). The symbol ``$\otimes$'' indicates the Kronecker product. Following Lange et al. (1989) we have, for the $q$-variate Student $t$ distribution with $\nu$ degrees of freedom,  $t_q({\mu}, {\Sigma}, \nu)$, that $W_g(u) = -(\nu + q )/\{2(\nu + u)\}$. Following G\'omez et al. (1998) we have, for the $q$-variate power exponential $PE_q({\mu}, \delta, \lambda)$ with shape parameter $\lambda>0$ and $u\neq0$, that $W_g(u) = - \lambda u^{\lambda - 1}/2$, $\lambda \neq 1/2$. In addition, we have

\begin{equation*} \label{matrix-Mi}
{M}_i =
\begin{bmatrix}
\frac{4\psi_{i(2,1)}}{q_{i}}{\Sigma}_i & {0}\\
{0} &   2c_i{\Sigma}_i\otimes{\Sigma}_i
\end{bmatrix}
+(c_i - 1)
\begin{bmatrix}
{0} & {0}\\
{0} &  \vvec({\Sigma}_i)\vvec({\Sigma}_i)^{\top}
\end{bmatrix},
\end{equation*} 
where $c_i = 4\psi_{i(2,2)}/\{q_{i}(q_{i} + 2)\}$, $\psi_{i(2,1)} = E(W_{g}^2(r_{i})r_{i})$ and $\psi_{i(2,2)} =$ \break $E(W_{g}^2(r_{i})r_{i}^{2})$,
with $r_{i} = ||{L}_{i}||^2$, ${L}_{i}\sim El_{q_{i}}({0},{I}_{q_{i}})$. Here, we assume that $g(u)$ is such that $\psi_{i(2,1)}$ and $\psi_{i(2,2)}$ 
exist and are finite for all $i=1, \ldots, n$. One can verify these results by using standard differentiation techniques and some standard matrix operations. 

The values of $\psi_{i(l,k)}$ are obtained from solving the following one-dimensional integrals (Lange et al. 1989):
\begin{equation}\label{psi}
\psi_{i(l,k)} = \int_{0}^{\infty} W_g(s^2)^l g(s^2)r^{q_i + 2k-1}c_{q_i} \mbox{d} s,
\end{equation} 
\noindent where $c_{q_i}=2\pi^{\frac{q_i}{2}}/\Gamma(\frac{q_i}{2})$ is the surface area of the unit sphere in $\mathbb{R}^{q_i}$ and $\Gamma(a)$ is the well-known gamma function. One can find these quantities for many distributions simply solving (\ref{psi}) algebraically or numerically. Table 2 shows these quantities for the normal, Cauchy, Student $t$ and power exponential distributions.

\begin{table}[!h]\label{Tab:psi}
\caption{Functions $\psi_{i(2,1)}$, $\psi_{i(2,2)}$, $\psi_{i(3,2)}$, $\psi_{i(3,3)}$ for normal, Cauchy, Student $t$ and power exponential (P.E.) distributions.}
\label{tab2}
{
\begin{tabular}{lcccc}
\hline\noalign{\smallskip}
       & $\psi_{i(2,1)}$ & $\psi_{i(2,2)}$           & $\psi_{i(3,2)}$           & \!\!\!\!\!\!\!\!\!\!\!\!\!\!\!\!\!\!\!\!\!\!\!\! \!\!\!\!\!\!\!\!\!\!\!\!\!\!\!\!\!\!\!\!\!\!\!\!\!\!\!\!\!\!\!\! $\psi_{i(3,3)}$  \\\hline
normal & $\frac{q_i}{4}$ & $\frac{q_i(q_i + 2)}{4}$  & $-\frac{q_i(q_i + 2)}{8}$ & $-\frac{q_i(q_i + 2)(q_i + 4)}{8}$ \ \ \ \ \ \ \ \ \ \ \ \ \ \ \ \ \ \ \ \ \ \ $q_i \geq 1$ \\ \\ \\
Cauchy &  $\frac{q_i(q_{i} + 1)}{4(q_{i} + 3)}$ &  $\frac{q_{i}(q_{i} + 2)(q_{i} + 1)}{4(q_{i} + 3)}$&   $-\frac{q_i(q_i+2)(q_i+1)^2}{8(q_i+3)(q_i+5)}$ &  $-\frac{q_i(q_i+2)(q_i+4)(q_i+1)^2 }{8(q_i+3)(q_i+5)}$ \ \ \ \ \ \ \ $q_i \geq 1$  \\  \\  \\
Student $t$ &  $\frac{q_i(q_{i} + \nu)}{4(q_{i} + \nu + 2)}$ &  $\frac{q_{i}(q_{i} + 2)(q_{i} + \nu)}{4(q_{i} + \nu + 2)}$&   $-\frac{q_i(q_i+2)(q_i+\nu)^2}{8(q_i+2+\nu)(q_i+4+\nu)}$ &  $-\frac{q_i(q_i+2)(q_i+4)(q_i+\nu )^2 }{8(q_i+2+\nu )(q_i+4+\nu)}$ \ \ \ \ \ \ \ $q_i \geq 1$ \\ \\ \\
%
P.E.& $\frac{\lambda^2\Gamma(\frac{4\lambda - 1}{2\lambda})}{2^{1/\lambda}\Gamma(\frac{1}{2\lambda})}$&   $\frac{2\lambda + 1}{4}$ & $- \frac{\lambda^3\Gamma(\frac{6\lambda - 1}{2\lambda})}{2^{1/\lambda}\Gamma(\frac{1}{2\lambda})}$&   $-\frac{(2\lambda + 1)(4\lambda + 1) }{8}$ \ \ \ \ \ \ \ \ \ \ $q_i  = 1$, \ $\lambda > \frac{1}{4}$ \\ \\
P.E.& $\frac{\lambda^2\Gamma(\frac{q_{i} - 2}{2\lambda} + 2)}{2^{1/\lambda}\Gamma(\frac{q_{i}}{2\lambda})}$&   $\frac{q_{i}(2\lambda + q_{i})}{4}$ & $- \frac{\lambda^3\Gamma(\frac{q_{i} - 2}{2\lambda} + 3)}{2^{1/\lambda}\Gamma(\frac{q_{i}}{2\lambda})}$&   $-\frac{q_i(2\lambda + q_i)(4\lambda +q_i) }{8}$ \ \  $q_i \geq 2$, \ $\lambda > 0$ \\
\noalign{\smallskip}\hline
\end{tabular}
}
\end{table}

It is important to remark that the $\psi_{i(l, k)}$'s may involve unknown quantities (for instance, the degrees of freedom $\nu$ of the Student $t$ distribution and the shape parameter $\lambda$ of the power exponential distribution). One may want to estimate these quantities via maximum likelihood estimation. Here, we consider these as known quantities for the purpose of keeping the robustness property of some distributions. Lucas (1997) shows that the protection against ``large'' observations is only valid when the degrees of freedom parameter is kept fixed for the Student $t$ distribution. Therefore, the issue of estimating these quantities is beyond of the main scope of this paper. In practice, one can use model selection procedures to choose the most appropriate values of such unknown parameters. 

Notice that, in the Fisher information matrix $K(\theta)$, the matrix ${M}$ carries all the information about the adopted distribution, while ${F}$ and ${H}$ contain the information about the adopted model. Also, $K(\theta)$ has a quadratic form that can be computed through simple matrix operations. Under the normal case, $v_{i} = 1$, ${M}={H}^{-1}$ and hence $\widetilde{H} = {H}$. 

The Fisher scoring method can be used to estimate ${\theta}$ by iteratively solving the equation 

\begin{equation} \label{Fisher-Scoring}
 ({F}^{(m)\top}\widetilde{H}^{(m)}{F}^{(m)}){\theta}^{(m+1)} =
  {F}^{(m)\top} \widetilde{H}^{(m)}{s}^{*(m)}, \quad  m = 0, 1,\ldots,
\end{equation}
where the quantities with the upper index ``$(m)$" are evaluated at $\widehat{{\theta}}$, $m$ is the iteration counter and $${s}^{*(m)} = {F}^{(m)}{\theta}^{(m)} + {H}^{-1(m)}{M}^{-1(m)}{s}^{(m)}.$$ \noindent Each loop, through the iterative scheme~(\ref{Fisher-Scoring}), consists of an iterative re-weighted least squares algorithm to optimize the log-likelihood~(\ref{log-likelihood}). Thus, (\ref{ScoreFisher2}) and (\ref{Fisher-Scoring}) agree with the corresponding equations derived in Patriota and Lemonte (2009). Observe that, despite the complexity and generality of the postulated model, expressions (\ref{ScoreFisher2}) and (\ref{Fisher-Scoring}) are very simple and friendly. 
 
Now, we can give the main result of the paper.
\vspace{0.2cm}
 
\noindent{\bf Theorem 3.1.} The second-order bias vector $B_{\widehat{\theta}}(\theta)$ under model (\ref{MainModel}) is given by
\begin{equation}\label{BIAS-vector}
B_{\widehat{\theta}}(\theta) = (F^{\top}\widetilde{{H}}F)^{-1} F^{\top}\widetilde{{H}}\xi,
\end{equation}
where ${\xi} = (\Phi_1,\ldots,\Phi_p)\vvec(({F}^{\top} \widetilde{H}F)^{-1})$, $\Phi_r = (\Phi_{1(r)}^{\top}, \ldots \Phi_{n(r)}^{\top})^{\top}$, and
$\Phi_{i(r)}$ is given in the Appendix. 

\vspace{0.5cm}
\noindent {\bf Proof:} See the Appendix.
\vspace{0.5cm}

In many models the location vector and the scale matrix do not have parameters in common, i.e., ${\mu}_i = {\mu}_i({\theta}_1,{x}_{i})$ and ${\Sigma}_i =$ ${\Sigma}_i({\theta}_2,{w}_{i})$, where ${\theta} = ({\theta}_1^{\top}, {\theta}_2^{\top} )^{\top}$. Therefore, $F = \mbox{block--diag}\{F_{\theta_1}, F_{\theta_2}\}$ and the parameter vectors ${\theta}_1$ and ${\theta}_2$ will be orthogonal (Cox and Reid 1987). This happens in mixed models, nonlinear models, among others. However, in errors-in-variables and factor analysis models orthogonality does not hold. Model (\ref{MainModel}) is general enough to encompass a large number of models even those that do not have orthogonal parameters.

\vspace{0.2cm}
 
\noindent{\bf Corollary 3.1.} When ${\mu}_i = {\mu}_i({\theta}_1,{x}_{i})$ and ${\Sigma}_i = {\Sigma}_i({\theta}_2,{w}_{i})$, where ${\theta} = ({\theta}_1^{\top}, {\theta}_2^{\top} )^{\top}$ the second-order bias vector of $\widehat{\theta}_1$ and $\widehat{\theta}_2$ are given by
\begin{equation*}\label{BIAS-vector1}
B_{\widehat{\theta}_1}(\theta) = (F_{\theta_1}^{\top} \widetilde{H}_{1} F_{\theta_1})^{-1} F_{\theta_1}^{\top} \widetilde{H}_{1} \xi_{1} 
\end{equation*}
and
\begin{equation*}\label{BIAS-vector2}
B_{\widehat{\theta}_2}(\theta) = (F_{\theta_2}^{\top} \widetilde{H}_{2} F_{\theta_2})^{-1} F_{\theta_2}^{\top} \widetilde{H}_{2} \xi_{2}, \end{equation*}
\noindent respectively. The quantities $F_{\theta_1}$, $F_{\theta_2}$, $\widetilde{H}_{1}$, $\widetilde{H}_{2}$, $\xi_{1}$ and $\xi_{2}$ are defined in the Appendix.

\vspace{0.5cm}
\noindent {\bf Proof:} See the Appendix.
\vspace{0.5cm}

\vspace{0.2cm}
\noindent Formula~(\ref{BIAS-vector}) says that, for any particular model of the general multivariate elliptical class of models (\ref{MainModel}), it is always possible to express the bias of $\widehat{{\theta}}$ as the solution of an weighted least-squares regression. Also, if $z_i \sim N_{q_i}(0, \Sigma_i)$ then $c_i = -\widetilde{\omega}_{i} = 1$, $\eta_{1i} = 0$, $\eta_{2i} = -2$, $\widetilde{H} = H$, 
\[J_{i(r)} = 
\begin{pmatrix}
0 \\ 
2(I_{q_i}\otimes a_{i(r)})D_i
\end{pmatrix},
\] and formula~(\ref{BIAS-vector}) reduces to the one obtained by Patriota and Lemonte (2009).

Theorem 3.1 implies that all one needs to compute bias-corrected and bias-reduced MLEs in the general elliptical model is: (i) the first and second derivatives of the location vector $\mu_i$ and the scale matrix $\Sigma_i$ with respect to all the parameters; (2) the derivatives $W_g(u)$; (3) some moments involving the chosen elliptical distribution (these moments are given in Table 2 for some elliptical distributions). With these quantities, the matrices in (\ref{BIAS-vector}) can be computed and the bias vector can be computed through an weighted least-squares regression.

\section{Special models}\label{special-cases}

In this section, we present four important particular cases of the main model (\ref{MainModel}). All special cases presented in Patriota and Lemonte (2009) are also special cases of the general multivariate elliptical model defined in this paper.

\subsection{Heteroscedastic nonlinear models}

Consider the univariate heteroscedastic nonlinear model defined by
\begin{equation*}\label{ModelNonLinHetero}
Y_i =  f({x}_i, {\alpha}) + e_i, \ \ i = 1, 2, \ldots, n,
\end{equation*}
\noindent where $Y_i$ is the response, ${x}_i$ is a column vector of explanatory variables, ${\alpha}$ is a column vector $p_1 \times 1$ of unknown parameters and $f$ is a nonlinear function of ${\alpha}$. Assume that $e_1, e_2, \ldots, e_n$ are independent, with $e_i \sim El(0, \sigma_i^2)$. Here $\sigma_i^2 = \sigma_i^2(\gamma) = h(\omega_i^\top \gamma)$, where $\gamma$ is a $p_2 \times 1$ vector of unknown parameters. Then

\begin{equation*}\label{DistHeteroNonLin}
{Y}_{i} \stackrel{ind}{\sim}El(f({x}_i, \alpha), \sigma_i^2),
\end{equation*}

\noindent which is a special case of (\ref{MainModel}) with $\theta = (\alpha^{\top}, \gamma^{\top})^{\top}$, ${\mu}_{i} = f({x}_i, \alpha)$ and ${\Sigma}_i = \sigma_i^2$. Here $El$ stands for $El_1$. Notice that for the heteroscedastic linear model $f({x}_i, {\alpha}) = {x}_i^{\top} {\alpha}$.

The second-order bias vector $B_{\widehat{\theta}}(\theta)$ comes from (\ref{BIAS-vector}), which depends on derivatives of $f({x}_i, \alpha)$ and $\sigma_i^2$ with respect to the parameter vector ${\theta}$. Also, it depends on the quantities $\psi_{i(2,1)}, \psi_{i(2,2)}, \psi_{i(3,2)}$, $\psi_{i(3,3)}$ (see Table 2) and $W_g(u_i)$ containing information about the adopted distribution.

\subsection{Nonlinear mixed-effects model}

One of the most important examples is the nonlinear mixed-effects model introduced by Lange et al. (1989) and studied under the assumption of a Student $t$ distribution. Let
\[Y_i =  \mu_i({x}_i, {\alpha}) + Z_i b_i + u_i,\]
where $Y_i$ is the $q_i \times 1$ vector response, $\mu_i$ is a $q_i$-dimensional nonlinear function of ${\alpha}$, ${x}_i$ is a vector of nonstochastic covariates, ${Z}_i$ is a matrix of known constants, ${\alpha}$ is a $p_1 \times 1$ vector of unknown parameters and ${b}_i$ is an $r \times 1$ vector of unobserved random regression coefficients. Assume that,
\begin{equation*}\label{Matrizbiui}
\left(\begin{array}{c} b_i \\ u_i  \\ \end{array}\right) \sim El_{r+q_i}\left(\left[\begin{array}{c} 0 \\ 0  \\ \end{array}\right], \left[\begin{array}{cc} \Sigma_{b}(\gamma_1) & 0  \\ 0 & R_i(\gamma_2) \\ \end{array}\right]\right),
\end{equation*}
where $\gamma_1$ is a $p_2$-dimensional vector of unknown parameters and $\gamma_2$ is a $p_3 \times 1$ vector of unknown parameters. Furthermore, the vectors $(b_1, u_1)^{\top}$, $(b_2, u_2)^{\top}$, $\ldots$, $(b_n, u_n)^{\top}$ are independent. Therefore, the marginal distribution of the observed vector is 

\begin{equation}\label{general-mixed}
Y_i \sim El_{q_i}\left(\mu_i({x}_i, {\alpha}); {\Sigma}_i({Z}_i, {\gamma})\right),
\end{equation} where $\gamma = (\gamma_1^{\top}, \gamma_2^{\top})^{\top}$ and ${\Sigma}_i(Z_i,{\gamma}) = 
{Z}_i \Sigma_b(\gamma_1){Z}_i^{\top} + R_i(\gamma_2)$. Equation (\ref{general-mixed}) is a special case of (\ref{MainModel}) with $\theta = (\alpha^{\top}, \gamma^{\top})^{\top}$, ${\mu}_{i} = \mu_i({x}_i, {\alpha})$ and ${\Sigma}_i = {\Sigma}_i({Z}_i, {\gamma})$. From (\ref{BIAS-vector}) one can compute the bias vector $B_{\widehat{\theta}}(\theta)$.

\subsection{Errors-in-variables model}

Consider the model

\begin{equation*}\label{error-model}
{x}_{1i} = {\beta}_0 + {\beta}_1 {x}_{2i} + {q}_i,\quad i=1,\ldots, n,
\end{equation*} 
where ${x}_{1i}$ is a $v\times 1$ latent response vector, ${x}_{2i}$ is a $m\times 1$ latent vector of covariates, ${\beta}_0$ is a $v\times 1$ vector of intercepts, ${\beta}_1$ is a $v \times m$ matrix of slopes, and ${q}_i$ is the equation error having a multivariate elliptical distribution with location vector zero and scale matrix ${\Sigma}_{{q}}$. The variables ${x}_{1i}$ and ${x}_{2i}$ are not directly observed, instead surrogate variables ${X}_{1i}$ and ${X}_{2i}$ are measured with the following additive structure: 
\begin{equation}\label{error-model-O}
{X}_{1i} = {x}_{1i} + {\delta}_{{x}_{1i}} \quad \mbox{and} \quad {X}_{2i} ={x}_{2i} + {\delta}_{{x}_{2i}}.
\end{equation}
The random quantities $x_{2i}$, $q_i$, ${\delta}_{{x}_{1i}}$ and ${\delta}_{{x}_{2i}}$ are assumed to follow an elliptical distribution given by
\[
\begin{pmatrix}
 x_{2i}\\
 q_i\\
 {\delta}_{{x}_{1i}}\\
 {\delta}_{{x}_{2i}}\\
\end{pmatrix}
\stackrel{ind}{\sim}El_{2v+2m}
\begin{bmatrix}
\begin{pmatrix}
\mu_{x_2}\\
0\\
 {0}\\
 {0}\\
\end{pmatrix},
\begin{pmatrix}
\Sigma_{x_2} & 0     & 0                 & 0\\
0        & \Sigma_q & 0                 &  0\\
0        &      0   & {\tau}_{{x}_{1i}} & {0} \\
0        &      {0} & 0                 & {\tau}_{{x}_{2i}}\\
\end{pmatrix}
\end{bmatrix},
\] where the matrices ${\tau}_{xi}$ and ${\tau}_{zi}$ are known for all $i=1, \ldots, n$. These ``known" matrices may be attained, for example, through an analytical treatment of the data collection mechanism, replications, machine precision, etc (Kulathinal et al. (2002)).

Therefore, the observed vector $Y_i =(X_{1i}^{\top}, X_{2i}^{\top})^{\top}$ has marginal distribution given by 
\begin{equation}\label{ModErros}
Y_i \stackrel{ind}{\sim}El_{v+m} (\mu(\theta), \Sigma_i(\theta))
\end{equation}
\noindent with
\[\mu(\theta) =
\begin{pmatrix}
\beta_0 + {\beta}_1 {\mu}_{{x_2}}\\
\mu_{x_2}
\end{pmatrix} \quad \mbox{and} \quad 
\Sigma_i(\theta)=
\begin{pmatrix}
\beta_1\Sigma_{x_2}\beta_1^{\top}+ \Sigma_q +
\tau_{{x}_{1i}} & {\beta}_1 \Sigma_{{x_{2}}}\\
{\Sigma}_{x_{2}}{\beta}_1^{\top} &{\Sigma}_{x_{2}} + \tau_{{x}_{2i}}
\end{pmatrix},
\] where ${\theta} = ({\beta}_0^{\top}, \vvec({\beta}_1)^{\top},{\mu}_{x_2}^{\top}, \mbox{vech}({\Sigma}_{x_2})^{\top}, \mbox{vech}(\Sigma_q)^{\top})^{\top}$, ``vech" operator transforms a symmetric matrix into a vector by stacking into columns its diagonal and superior diagonal elements. The mean vector $(\theta)$ and the covariance-variance matrix $\Sigma_i(\theta)$ of observed variables have the matrix ${\beta}_1$ in common, i.e., they share $mv$ parameters. Kulathinal et al. (2002) study the linear univariate case ($v=1$, $m=1$).

Equation (\ref{ModErros}) is a special case of (\ref{MainModel}) with $q_i = v + m$, $\theta = (\alpha^{\top}, \gamma^{\top})^{\top}$, ${\mu}_{i} = \mu_i(\theta)$ and ${\Sigma}_i = {\Sigma}_i(\theta)$. In this case, a programming language or software that can perform operations on vectors and matrices, e.g. {\tt Ox} (Doornik, 2013) and {\tt R} (Ihaka and Gentleman, 1996), can be used to obtain the bias vector $B_{\widehat{\theta}}(\theta)$ from (\ref{BIAS-vector}).

\subsection{Log-symmetric regression models}

Let $T$ be a continuous positive random variable with probability density function

\begin{equation}\label{ModelLogSym1}
f_T(t; \eta, \phi, g) = \frac{1}{\sqrt{\phi} t} g\left(\log^2 \left[\left(\frac{t}{\eta}\right)^{\frac{1}{\sqrt{\phi}}}\right]\right), \ \eta > 0, \ \phi > 0,
\end{equation}
where $g$ is the density generating function of a univariate elliptical distribution, and we write $T \sim LS(\eta, \phi, g)$. Vanegas and Paula (2014) called the class of distribution in (\ref{ModelLogSym1}) the log-symmetric class of distributions. It includes log-normal, log-Student $t$, log-power-exponential distributions, among many others, as special cases. It is easy to verify that $\log(T)$ has a univariate elliptical distribution (i.e., symmetric distribution) with location parameter $\mu = \log(\eta)$ and scale parameter $\phi$. The parameter $\eta$ is the median of $T$, and $\phi$ can be interpreted as a skewness or relative dispersion parameter. 

Vanegas and Paula (2015) defined and studied semi-parametric regression models for a set $T_1, T_2, \ldots, T_n$ with $T_i \sim LS(\eta_i, \phi_i, g)$ with $\eta_i > 0$ and $\phi_i > 0$ following semi-parametric regression structures. Here we assume parametric specification for $\eta_i$ and $\phi_i$ as $\eta_i = \eta_i(x_i, \alpha)$ and $\phi_i = \phi_i(\omega_i, \gamma)$.

Hence,

\begin{equation}\label{ModelLogSym2}
Y_i = \log(T_i) \stackrel{ind}{\sim}El\left(\mu_i(x_i, \alpha), \phi_i(\omega_i, \gamma)\right),
\end{equation}
where $\mu_i(x_i, \alpha) = \log(\eta(x_i, \alpha))$. Therefore, (\ref{ModelLogSym2}) is a special case of the general elliptical model (\ref{MainModel}), and formula (\ref{BIAS-vector}) applies.

\section{Simulation results}\label{simulation}

In this section, we shall present the results of Monte Carlo simulation experiments in which we evaluate the finite sample performances of the original MLEs and their bias-corrected and bias-reduced versions. The simulations are based on the univariate nonlinear model without random effects (Section 4.1) and the errors-in-variables model presented in Section 4.2, when $Y_i$ follows a normal distribution, a Student $t$ distribution with $\nu$ degrees of freedom, or a power exponential distribution with shape parameter $\lambda$. For all the simulations, the number of Monte Carlo replications is 10,000 (ten thousand) and they have been performed using the \texttt{Ox} matrix programming language (Doornik, 2013).

First consider the model described in (\ref{general-mixed}) with $q_i=1$, $Z_i = 0$, $\Sigma_i = \sigma^2$ and
\begin{equation}\label{nonlinear-model}
\mu_i({\alpha}) = \mu_i({x}_i, {{\alpha}}) = \alpha_1 + \frac{\alpha_2}{1 + \alpha_3 x_i^{\alpha_4}},\quad i = 1,\ldots, n.
\end{equation}

\noindent Here the unknown parameter vector is ${\theta} = (\alpha_1, \alpha_2, \alpha_3, \alpha_4, \sigma^2)^\top$. The values of $x_i$ were obtained as random draws from the uniform distribution $U(0,100)$. The sample sizes considered are $n = 10, 20, 30, 40$ and $50$. The parameter values are $\alpha_1 = 50$, $\alpha_2 = 500$, $\alpha_3 = 0.50$, $\alpha_4 = 2$ and $\sigma_i^2 = 200$. For the Student $t$ distribution, we fixed the degrees of freedom at $\nu = 4$, and for the power exponential model the shape parameter is fixed at $\lambda = 0.8$.

Tables 3-4 present the bias, and the root mean squared errors ($\sqrt{MSE}$) of the maximum likelihood estimates, the bias-corrected estimates and the bias-reduced estimates for the nonlinear model with normal and Student $t$ distributed errors, respectively. To save space, the corresponding results for the power exponential model are not shown.\footnote{All the omitted tables in this paper are presented in a supplement available from the authors upon request.} We note that the bias-corrected estimates and the bias-reduced estimates are less biased than the original MLE for all the sample sizes considered. For instance, when $n=20$ and the errors follow a Student $t$ distribution (see Table 4) the estimated biases of ${\widehat \sigma}^2$ are $-41.24$ (MLE), $-12.30$ (bias-corrected) and $-4.55$ (bias-reduced). For the normal case with $n = 10$ (see Table 3), the estimated biases of ${\widehat \alpha}_2$ are $2.16$ (MLE), $0.70$ (bias-corrected) and $-0.27$ (bias-reduced). We also observe that the bias-reduced estimates are less biased than the bias-corrected estimates in most cases. As $n$ increases, the bias and the root mean squared error of all the estimators decrease, as expected. Additionally, we note that the MLE of ${\alpha}_2$ has $\sqrt{MSE}$ larger than those of the modified versions. For the estimation of $\sigma^2$, $\sqrt{MSE}$ is smaller for the original MLE. In other cases, we note that the estimators have similar root mean squared errors.

\begin{table}[!htp]
{\caption{Biases and $\sqrt{MSE}$ of the maximum likelihood estimate and its adjusted versions; nonlinear model; normal distribution.}} 
\label{tab3}
\begin{center}
\begin{tabular}{rrrrrrrrrrrr}
\hline\noalign{\smallskip}
   &                 && \multicolumn{2}{c}{MLE} && \multicolumn{2}{c}{Bias-corrected MLE} && \multicolumn{2}{c}{Bias-reduced MLE}  \\ 
\noalign{\smallskip}   
 \cline{4-5}\cline{7-8}\cline{10-11}
\noalign{\smallskip}
$n$&${\theta}$ &&    Bias  &$\sqrt{MSE}$ &&    Bias  &$\sqrt{MSE}$  &&     Bias  &$\sqrt{MSE}$  \\ 
\noalign{\smallskip}
\cline{1-2} \cline{4-5} \cline{7-8}\cline{10-11}
\noalign{\smallskip}
   & $\alpha_1$      && \    $-0.29$ & \ \      $6.69$ && \    $-0.13$ & \ \      $6.67$ && \   $-0.01$ & \ \      $6.67$ \\       
   & $\alpha_2$      && \ \   $2.16$ & \!      $20.07$ && \ \   $0.70$ & \!      $19.40$ && \   $-0.27$ & \ \     $19.06$ \\       
10 & $\alpha_3$      && \ \   $0.01$ & \ \      $0.13$ && \ \   $0.00$ & \ \      $0.12$ && \ \  $0.00$ & \ \      $0.12$ \\       
   & $\alpha_4$      && \ \   $0.03$ & \ \      $0.30$ && \ \   $0.01$ & \ \      $0.29$ && \   $-0.00$ & \ \      $0.29$ \\       
   & $\sigma^2$      &&     $-80.05$ & \!\!\! $106.44$ &&     $-32.06$ & \!\!\! $103.32$ && \ \  $9.09$ & \!\!\! $128.72$ \\ 
   \noalign{\smallskip} \hline \noalign{\smallskip}
   & $\alpha_1$      && \    $-0.08$ & \ \      $4.07$ && \   $-0.01$  & \ \      $4.07$ && \ \  $0.01$ & \ \      $4.07$ \\       
   & $\alpha_2$      && \ \   $0.66$ & \!      $17.94$ && \   $-0.08$  & \!      $17.84$ && \   $-0.27$ & \!      $17.82$ \\       
20 & $\alpha_3$      && \ \   $0.00$ & \ \      $0.09$ && \ \   $0.00$ & \ \      $0.09$ && \   $-0.00$ & \ \      $0.09$ \\       
   & $\alpha_4$      && \ \   $0.02$ & \ \      $0.21$ && \ \   $0.01$ & \ \      $0.20$ && \ \  $0.00$ & \ \      $0.20$ \\       
   & $\sigma^2$      &&     $-40.07$ & \!      $69.73$ && \    $-8.09$ & \!      $68.95$ && \ \  $0.86$ & \!      $72.02$ \\ 
   \noalign{\smallskip} \hline \noalign{\smallskip}
   & $\alpha_1$      && \    $-0.10$ & \ \      $3.11$ && \    $-0.04$ & \ \      $3.10$ && \   $-0.02$ & \ \      $3.10$ \\       
   & $\alpha_2$      && \ \   $0.71$ & \!      $17.24$ && \    $-0.05$ & \!      $17.15$ && \   $-0.18$ & \!      $17.13$ \\       
30 & $\alpha_3$      && \ \   $0.00$ & \ \      $0.09$ && \    $-0.00$ & \ \      $0.09$ && \   $-0.00$ & \ \      $0.09$ \\       
   & $\alpha_4$      && \ \   $0.02$ & \ \      $0.20$ && \ \   $0.00$ & \ \      $0.19$ && \ \  $0.00$ & \ \      $0.19$ \\       
   & $\sigma^2$      &&     $-26.41$ & \!      $55.26$ && \    $-3.26$ & \!      $55.11$ && \ \  $0.82$ & \!      $56.32$ \\ 
   \noalign{\smallskip} \hline \noalign{\smallskip}
   & $\alpha_1$      && \    $-0.08$ & \ \      $2.69$ && \    $-0.02$ & \ \      $2.69$ && \   $-0.01$ & \ \      $2.69$ \\       
   & $\alpha_2$      && \ \   $0.83$ & \!      $16.80$ && \ \   $0.09$ & \!      $16.70$ && \ \  $0.01$ & \!      $16.69$ \\       
40 & $\alpha_3$      && \ \   $0.00$ & \ \      $0.09$ && \ \   $0.00$ & \ \      $0.09$ && \ \  $0.00$ & \ \      $0.09$ \\       
   & $\alpha_4$      && \ \   $0.02$ & \ \      $0.19$ && \ \   $0.00$ & \ \      $0.18$ && \   $-0.00$ & \ \      $0.18$ \\       
   & $\sigma^2$      &&     $-20.04$ & \!      $47.26$ && \    $-2.04$ & \!      $47.13$ && \ \  $0.33$ & \!      $47.74$ \\ 
   \noalign{\smallskip} \hline \noalign{\smallskip}
   & $\alpha_1$      && \    $-0.08$ & \ \      $2.39$ && \    $-0.03$ & \ \      $2.38$ && \   $-0.02$ & \ \      $2.38$ \\       
   & $\alpha_2$      && \ \   $1.07$ & \!      $14.25$ && \ \   $0.30$ & \!      $14.12$ && \ \  $0.23$ & \!      $14.11$ \\       
50 & $\alpha_3$      && \ \   $0.00$ & \ \      $0.08$ && \ \   $0.00$ & \ \      $0.08$ && \ \  $0.00$ & \ \      $0.08$ \\       
   & $\alpha_4$      && \ \   $0.01$ & \ \      $0.19$ && \ \   $0.00$ & \ \      $0.18$ && \   $-0.00$ & \ \      $0.18$ \\       
   & $\sigma^2$      &&     $-15.93$ & \!      $41.41$ && \    $-1.21$ & \!      $41.30$ && \ \  $0.36$ & \!      $41.67$ \\ 
\noalign{\smallskip}\hline
\end{tabular}	
\end{center}
\end{table}

\begin{table}[!htp]
{\caption{Biases and $\sqrt{MSE}$ of the maximum likelihood estimate and its adjusted versions; nonlinear model; Student $t$ distribution.}} 
\label{tab4}
\begin{center}
\begin{tabular}{rrrrrrrrrrrr}
\hline\noalign{\smallskip}
   &                 && \multicolumn{2}{c}{MLE} && \multicolumn{2}{c}{Bias-corrected MLE} && \multicolumn{2}{c}{Bias-reduced MLE}  \\ 
\noalign{\smallskip}   
 \cline{4-5}\cline{7-8}\cline{10-11}
\noalign{\smallskip}
$n$&${\theta}$ &&    Bias  &$\sqrt{MSE}$ &&    Bias  &$\sqrt{MSE}$  &&     Bias  &$\sqrt{MSE}$  \\ 
\noalign{\smallskip}
\cline{1-2} \cline{4-5} \cline{7-8}\cline{10-11}
\noalign{\smallskip}
   & $\alpha_1$      && \     $-0.51$ & \ \      $8.66$ && \    $-0.31$ & \ \      $8.63$ && \   $-0.20$ & \ \      $8.56$ \\       
   & $\alpha_2$      && \ \    $3.34$ & \!      $28.47$ && \ \   $1.39$ & \!      $27.34$ && \ \  $2.05$ & \!      $27.67$ \\       
10 & $\alpha_3$      && \ \    $0.01$ & \ \      $0.17$ && \ \   $0.00$ & \ \      $0.16$ && \ \  $0.01$ & \ \      $0.16$ \\       
   & $\alpha_4$      && \ \    $0.06$ & \ \      $0.42$ && \ \   $0.03$ & \ \      $0.39$ && \   $-0.01$ & \ \      $0.38$ \\       
   & $\sigma^2$      &&      $-93.18$ & \!\!\! $127.60$ &&     $-54.24$ & \!\!\! $130.73$ &&    $-17.40$ & \!\!\! $170.35$ \\ 
   \noalign{\smallskip}\hline\noalign{\smallskip}
   & $\alpha_1$      && \     $-0.17$ & \ \      $5.03$ && \!   $-0.07$ & \ \      $5.02$ && \   $-0.04$ & \ \      $5.01$ \\       
   & $\alpha_2$      && \ \    $2.01$ & \!      $25.64$ && \ \   $0.91$ & \!      $25.11$ && \ \  $1.29$ & \!      $24.98$ \\       
20 & $\alpha_3$      && \ \    $0.01$ & \ \      $0.14$ && \ \   $0.01$ & \ \      $0.14$ && \ \  $0.01$ & \ \      $0.13$ \\       
   & $\alpha_4$      && \ \    $0.04$ & \ \      $0.29$ && \ \   $0.01$ & \ \      $0.28$ && \ \  $0.00$ & \ \      $0.27$ \\       
   & $\sigma^2$      &&      $-41.24$ & \       $85.51$ &&     $-12.30$ & \       $89.41$ && \   $-4.55$ & \!      $93.08$ \\ 
   \noalign{\smallskip}\hline\noalign{\smallskip}
   & $\alpha_1$      && \     $-0.10$ & \ \      $3.81$ && \    $-0.01$ & \ \      $3.80$ && \ \  $0.01$ & \ \      $3.82$ \\       
   & $\alpha_2$      && \ \    $2.25$ & \!      $25.75$ && \ \   $1.13$ & \!      $25.34$ && \ \  $1.61$ & \!      $25.41$ \\       
30 & $\alpha_3$      && \ \    $0.01$ & \ \      $0.14$ && \ \   $0.01$ & \ \      $0.14$ && \ \  $0.01$ & \ \      $0.14$ \\       
   & $\alpha_4$      && \ \    $0.04$ & \ \      $0.29$ && \ \   $0.01$ & \ \      $0.27$ && \ \  $0.00$ & \ \      $0.26$ \\       
   & $\sigma^2$      &&      $-27.15$ & \!      $70.02$ && \    $-6.15$ & \!      $72.64$ && \   $-1.78$ & \!\!\! $107.53$ \\ 
   \noalign{\smallskip}\hline\noalign{\smallskip}
   & $\alpha_1$      && \     $-0.10$ & \ \      $3.27$ && \    $-0.02$ & \ \      $3.26$ && \   $-0.01$ & \ \      $3.26$ \\       
   & $\alpha_2$      && \ \    $1.82$ &         $24.94$ && \ \   $0.75$ &         $24.67$ && \ \  $1.18$ &         $24.78$ \\       
40 & $\alpha_3$      && \ \    $0.01$ & \ \      $0.12$ && \ \   $0.00$ & \ \      $0.12$ && \ \  $0.01$ & \ \      $0.12$ \\       
   & $\alpha_4$      && \ \    $0.03$ & \ \      $0.26$ && \ \   $0.01$ & \ \      $0.25$ && \ \  $0.00$ & \ \      $0.25$ \\       
   & $\sigma^2$      &&      $-20.38$ &         $60.43$ && \    $-4.01$ &         $62.21$ && \   $-1.82$ &         $62.98$ \\ 
   \noalign{\smallskip}\hline\noalign{\smallskip}
   & $\alpha_1$      && \     $-0.13$ & \ \      $2.86$ && \    $-0.05$ & \ \      $2.85$ && \   $-0.03$ & \ \      $2.85$ \\       
   & $\alpha_2$      && \ \    $1.48$ &         $18.86$ && \ \   $0.38$ &         $18.59$ && \ \  $0.24$ &         $18.46$ \\       
50 & $\alpha_3$      && \ \    $0.01$ & \ \      $0.11$ && \ \   $0.00$ & \ \      $0.11$ && \ \  $0.00$ & \ \      $0.11$ \\       
   & $\alpha_4$      && \ \    $0.02$ & \ \      $0.24$ && \ \   $0.00$ & \ \      $0.23$ && \ \  $0.00$ & \ \      $0.23$ \\       
   & $\sigma^2$      &&      $-15.40$ &         $53.99$ && \    $-1.94$ &         $55.56$ && \   $-0.43$ &         $56.11$ \\ 
   \noalign{\smallskip}\hline
\end{tabular}	 
\end{center}          			  
\end{table}

We now consider the errors-in-variables model described in (\ref{error-model-O}). The sample sizes considered are $n = 15, 25, 35$ and $50$. The parameter values are ${\beta}_0 = 0.70 \ {\bfsym 1}_{v \times 1}$, ${\beta}_1 = 0.40 \ {\bfsym 1}_{v \times m}$, ${\mu}_{x_2} = 70 \ {\bfsym 1}_{m \times 1}$, ${\Sigma}_{q} = 40 \ {\bfsym I}_{v \times 1}$ and ${\Sigma}_{x_2} = 250 \ {\bfsym I}_{m \times 1}$. Here, ${\bfsym 1}_{r \times s}$ is as $r \times s$ matrix of ones and $ {\bfsym I}_{r \times s}$ is the identity matrix with dimension $r \times s$. For the Student $t$ distribution, we fixed the degrees of freedom at $\nu = 4$ and, for power exponential model, the shape parameter was fixed at $\lambda = 0.7$. We consider $v \in \{1, 2\}$ and $m = 1$.

In Tables 5-6, we present the MLE, the bias-corrected estimates, the bias-reduced estimates, and corresponding estimated root mean squared errors for the Student $t$ and power exponential distributions, for the errors-in-variables model. The results for the normal distribution are not shown to save space. We observe that, in absolute value, the biases of the bias-corrected estimates and bias-reduced estimates are smaller than those of the original MLE for different sample sizes. Furthermore, the bias-reduced estimates are less biased than the bias-corrected estimates in most cases. This can be seen e.g. in Table 6 when $v = 1$, $m = 1$, $Y_i$ follows a power exponential distribution and $n=15$. In this case, the bias of the MLE, the bias-corrected estimate and the bias-reduced estimate of ${\Sigma}_q$ are $-4.92$, $-0.66$ and $-0.17$, respectively. When $Y_i$ follows a Student $t$ distribution, $n = 15$, $v = 1$ and $m = 1$ we observe the following biases of the estimates of ${\Sigma}_{x_2}$: $5.18$ (MLE), $2.91$ (bias-corrected) and $2.66$ (bias-reduced); see Table 5. We note that the root mean squared errors decrease with $n$. 

For the sake of saving space, the simulation results for the normal, Student $t$ and power exponential errors-in-variable models with $v=2$ and $m=1$ are not presented. Overall, our findings are similar to those reached for the other models.

\begin{table}[!htp]
{\caption{Biases and $\sqrt{MSE}$ of the maximum likelihood estimate and its adjusted versions; errors-in-variables model; $v = 1$ and $m = 1$; Student $t$ distribution.}} 
\label{tab7}
\begin{center}
\begin{tabular}{rrrrrrrrrrrr} 
\hline\noalign{\smallskip}
   &                && \multicolumn{2}{c}{MLE} && \multicolumn{2}{c}{Bias-corrected MLE} && \multicolumn{2}{c}{Bias-reduced MLE}  \\ 
\noalign{\smallskip}\cline{4-5}\cline{7-8}\cline{10-11}\noalign{\smallskip}
$n$&${\theta}$&&    Bias      &$\sqrt{MSE}$     &&    Bias   &$\sqrt{MSE}$    &&     Bias  &$\sqrt{MSE}$    \\ 
\noalign{\smallskip}\cline{1-2} \cline{4-5} \cline{7-8}\cline{10-11}\noalign{\smallskip}
   & $\beta_0$      && \      $-0.00$ & \ \       $9.90$ && \ \  $0.01$ & \ \      $9.90$ && \ \  $0.01$ & \ \      $9.89$ \\       
   & $\beta_1$      && \ \     $0.00$ & \ \       $0.14$ && \ \  $0.00$ & \ \      $0.14$ && \ \  $0.00$ & \ \      $0.14$ \\       
15 & $\mu_{x_2}$    && \ \     $0.05$ & \ \       $4.82$ && \ \  $0.05$ & \ \      $4.82$ && \ \  $0.05$ & \ \      $4.82$ \\       
   & $\Sigma_{x_2}$ && \ \     $5.18$ & \!\!\!  $129.34$ && \ \  $2.91$ & \!\!\! $128.12$ && \ \  $2.66$ & \!\!\! $127.86$ \\       
   & $\Sigma_q$     && \      $-3.64$ & \!       $19.52$ && \   $-0.68$ & \!      $20.72$ && \   $-0.42$ & \!      $20.85$ \\ 
   \noalign{\smallskip}\hline\noalign{\smallskip}
   & $\beta_0$      && \      $-0.02$ & \ \       $7.14$ && \   $-0.01$ & \ \      $7.14$ && \   $-0.01$ & \ \      $7.14$ \\       
   & $\beta_1$      && \ \     $0.00$ & \ \       $0.10$ && \ \  $0.00$ & \ \      $0.10$ && \ \  $0.00$ & \ \      $0.10$ \\       
25 & $\mu_{x_2}$    && \ \     $0.03$ & \ \       $3.69$ && \ \  $0.03$ & \ \      $3.69$ && \ \  $0.03$ & \ \      $3.69$ \\       
   & $\Sigma_{x_2}$ && \ \     $3.61$ & \!       $97.32$ && \ \  $2.25$ & \!      $96.76$ && \ \  $2.17$ & \!      $96.69$ \\       
   & $\Sigma_q$     && \      $-2.31$ & \!       $14.87$ && \   $-0.47$ & \!      $15.41$ && \   $-0.38$ & \!      $15.44$ \\ 
   \noalign{\smallskip}\hline\noalign{\smallskip}
   & $\beta_0$      && \      $-0.02$ & \ \       $5.93$ && \   $-0.02$ & \ \      $5.93$ && \   $-0.02$ & \ \      $5.93$ \\       
   & $\beta_1$      && \ \     $0.00$ & \ \       $0.08$ && \ \  $0.00$ & \ \      $0.08$ && \ \  $0.00$ & \ \      $0.08$ \\       
35 & $\mu_{x_2}$    && \      $-0.01$ & \ \       $3.12$ && \   $-0.01$ & \ \      $3.12$ && \   $-0.01$ & \ \      $3.12$ \\       
   & $\Sigma_{x_2}$ && \ \     $1.94$ & \!       $79.78$ && \ \  $0.98$ & \!      $79.45$ && \ \  $0.94$ & \!      $79.44$ \\       
   & $\Sigma_q$     && \      $-1.65$ & \!       $12.63$ && \   $-0.31$ & \!      $12.96$ && \   $-0.26$ & \!      $12.97$ \\ 
   \noalign{\smallskip}\hline\noalign{\smallskip}
   & $\beta_0$      && \      $-0.01$ & \ \       $4.92$ && \   $-0.01$ & \ \      $4.92$ && \   $-0.01$ & \ \      $4.92$ \\       
   & $\beta_1$      && \ \     $0.00$ & \ \       $0.07$ && \ \  $0.00$ & \ \      $0.07$ && \ \  $0.00$ & \ \      $0.07$ \\       
50 & $\mu_{x_2}$    && \ \     $0.01$ & \ \       $2.59$ && \ \  $0.01$ & \ \      $2.59$ && \ \  $0.01$ & \ \      $2.59$ \\       
   & $\Sigma_{x_2}$ && \ \     $1.04$ & \!       $65.50$ && \ \  $0.37$ & \!      $65.33$ && \ \  $0.36$ & \!      $65.33$ \\       
   & $\Sigma_q$     && \      $-1.18$ & \!       $10.53$ && \   $-0.24$ & \!      $10.78$ && \   $-0.21$ & \!      $10.79$ \\ 
   \noalign{\smallskip}\hline\noalign{\smallskip}
   
\end{tabular}	
\end{center}				           			  
\end{table}

\begin{table}[!htp]
{\caption{Biases and $\sqrt{MSE}$ of the maximum likelihood estimate and its adjusted versions; errors-in-variables model; $v = 1$ and $m = 1$; power exponential distribution.}} 
\label{tab8}
\begin{center}
\begin{tabular}{rrrrrrrrrrrr}
\hline\noalign{\smallskip}
   &                 && \multicolumn{2}{c}{MLE} && \multicolumn{2}{c}{Bias-corrected MLE} && \multicolumn{2}{c}{Bias-reduced MLE}  \\ 
\noalign{\smallskip}   
 \cline{4-5}\cline{7-8}\cline{10-11}
\noalign{\smallskip}
$n$&${\theta}$ &&    Bias  &$\sqrt{MSE}$ &&    Bias  &$\sqrt{MSE}$  &&     Bias  &$\sqrt{MSE}$  \\ 
\noalign{\smallskip}
\cline{1-2} \cline{4-5} \cline{7-8}\cline{10-11}
\noalign{\smallskip}
   & $\beta_0$      && \      $-0.12$ & \ \      $9.25$ && \   $-0.11$ & \ \      $9.25$ && \   $-0.11$ & \ \     $9.24$ \\       
   & $\beta_1$      && \ \     $0.00$ & \ \      $0.13$ && \ \  $0.00$ & \ \      $0.13$ && \ \  $0.00$ & \ \     $0.13$ \\       
15 & $\mu_{x_2}$    && \      $-0.02$ & \ \      $6.47$ && \   $-0.02$ & \ \      $6.47$ && \   $-0.02$ & \ \     $6.47$ \\       
   & $\Sigma_{x_2}$ && \      $-9.27$ & \!\!\! $103.32$ && \ \  $0.52$ & \!\!\! $107.51$ && \ \  $0.82$ & \!\!\!$107.64$ \\       
   & $\Sigma_q$     && \      $-4.92$ & \!      $15.67$ && \   $-0.66$ & \!      $17.55$ && \   $-0.17$ & \!     $17.76$ \\ 
   \noalign{\smallskip}\hline\noalign{\smallskip}
   & $\beta_0$      && \ \     $0.02$ & \ \      $6.83$ && \ \  $0.03$ & \ \      $6.83$ && \ \  $0.03$ & \ \     $6.83$ \\       
   & $\beta_1$      && \ \     $0.00$ & \ \      $0.09$ && \   $-0.00$ & \ \      $0.09$ && \   $-0.00$ & \ \     $0.09$ \\       
25 & $\mu_{x_2}$    && \      $-0.02$ & \ \      $4.98$ && \   $-0.02$ & \ \      $4.98$ && \   $-0.02$ & \ \     $4.98$ \\       
   & $\Sigma_{x_2}$ && \      $-5.60$ & \!      $80.20$ && \ \  $0.36$ & \!      $81.95$ && \ \  $0.47$ & \!     $81.99$ \\       
   & $\Sigma_q$     && \      $-3.04$ & \!      $12.94$ && \   $-0.36$ & \!      $13.49$ && \   $-0.18$ & \!     $13.54$ \\ 
   \noalign{\smallskip}\hline\noalign{\smallskip}
   & $\beta_0$      && \ \     $0.01$ & \ \      $5.59$ && \ \  $0.02$ & \ \      $5.58$ && \ \  $0.02$ & \ \     $5.58$ \\       
   & $\beta_1$      && \      $-0.00$ & \ \      $0.08$ && \   $-0.00$ & \ \      $0.08$ && \   $-0.00$ & \ \     $0.08$ \\       
35 & $\mu_{x_2}$    && \      $-0.04$ & \ \      $4.21$ && \   $-0.04$ & \ \      $4.21$ && \   $-0.04$ & \ \     $4.21$ \\       
   & $\Sigma_{x_2}$ && \      $-3.53$ & \!      $68.01$ && \ \  $0.77$ & \!      $69.10$ && \ \  $0.82$ & \!     $69.12$ \\       
   & $\Sigma_q$     && \      $-2.14$ & \!      $11.11$ && \   $-0.18$ & \!      $11.46$ && \   $-0.08$ & \!     $11.49$ \\ 
   \noalign{\smallskip}\hline\noalign{\smallskip}
   & $\beta_0$      && \ \     $0.03$ & \ \      $4.67$ && \ \  $0.03$ & \ \      $4.67$ && \ \  $0.03$ & \ \     $4.67$ \\       
   & $\beta_1$      && \      $-0.00$ & \ \      $0.06$ && \   $-0.00$ & \ \      $0.06$ && \   $-0.00$ & \ \     $0.06$ \\       
50 & $\mu_{x_2}$    && \      $-0.03$ & \ \      $3.52$ && \   $-0.03$ & \ \      $3.52$ && \   $-0.03$ & \ \     $3.52$ \\       
   & $\Sigma_{x_2}$ && \      $-2.83$ & \!      $56.89$ && \ \  $0.18$ & \!      $57.51$ && \ \  $0.21$ & \!     $57.52$ \\       
   & $\Sigma_q$     && \      $-1.51$ & \ \      $9.21$ && \   $-0.12$ & \ \      $9.41$ && \   $-0.07$ & \ \     $9.42$ \\ 
   \noalign{\smallskip}\hline  
\end{tabular}			
\end{center}				           			  
\end{table}

\section{Applications}\label{applications}

\subsection{Radioimmunoassay data}

Tiede and Pagano (1979) present a dataset, referred here as the radioimmunoassay data, obtained from the Nuclear Medicine Department at the Veterans Administration Hospital, Buffalo, New York. Lemonte and Patriota (2011) analyzed the data to illustrate the applicability of the elliptical models with general parameterization. Following Tiede and Pagano (1979) we shall consider the nonlinear regression model (\ref{nonlinear-model}), with $n = 14.$ The response variable is the observed radioactivity (count in thousands), the covariate corresponds to the thyrotropin dose (measured in micro-international units per milliliter) and the errors follow a normal distribution or a Student $t$ distribution with $\nu = 4$ degrees of freedom. We assume that the scale parameter is unknown for both models. In Table 7 we present the maximum likelihood estimates, the bias-corrected estimates, the bias-reduced estimates, and the corresponding estimated standard errors are given in parentheses. 
We note that all the estimates present smaller standard errors under the Student $t$ model than under the normal model (Table 7).

For all parameters, the original MLEs are very close to the bias-corrected MLE and the bias-reduced MLE when the Student $t$ model is used. However, under the normal model, significant differences in the estimates of $\alpha_1$ are noted. The estimates for $\alpha_1$ are $0.44$ (MLE), $0.65$ (bias-corrected MLE) and $1.03$ (bias-reduced MLE). 

\begin{table}[!htp]
\caption{Estimates and standard errors (given in parentheses); radioimmunoassay data.}
\label{tab12} 
\begin{center}
\begin{tabular}{cccc} 
\hline\noalign{\smallskip}
\multicolumn{4}{c}{Normal distribution}  \\ 
\noalign{\smallskip}\hline\noalign{\smallskip}   
${\theta}$        &       MLE         & Bias-corrected MLE   & Bias-reduced MLE    \\ 
\noalign{\smallskip}\hline\noalign{\smallskip}
  $\alpha_1$      & $0.44$ \ $(0.80)$ & $0.65$ \ $(0.99)$ & $1.03$ \ $(1.06)$ \\                               
  $\alpha_2$      & $7.55$ \ $(0.95)$ & $7.34$ \ $(1.16)$ & $6.91$ \ $(1.25)$ \\                               
  $\alpha_3$      & $0.13$ \ $(0.06)$ & $0.13$ \ $(0.06)$ & $0.13$ \ $(0.08)$ \\                               
  $\alpha_4$      & $0.96$ \ $(0.24)$ & $0.93$ \ $(0.28)$ & $0.95$ \ $(0.34)$ \\                               
  $\sigma^2$      & $0.31$ \ $(0.12)$ & $0.40$ \ $(0.15)$ & $0.50$ \ $(0.19)$ \\
  \noalign{\smallskip}\hline\noalign{\smallskip}
\multicolumn{4}{c}{Student $t$ distribution}  \\ 
\noalign{\smallskip}\hline\noalign{\smallskip}
${\theta}$  &       MLE         & Bias-corrected MLE   & Bias-reduced MLE    \\ 
\noalign{\smallskip}\hline\noalign{\smallskip}
  $\alpha_1$      & $0.90$ \ $(0.12)$ & $0.91$ \ $(0.13)$ & $0.90$ \ $(0.15)$ \\                               
  $\alpha_2$      & $7.09$ \ $(0.17)$ & $7.08$ \ $(0.19)$ & $7.07$ \ $(0.22)$ \\                               
  $\alpha_3$      & $0.09$ \ $(0.01)$ & $0.09$ \ $(0.01)$ & $0.09$ \ $(0.02)$ \\                               
  $\alpha_4$      & $1.31$ \ $(0.08)$ & $1.31$ \ $(0.09)$ & $1.29$ \ $(0.10)$ \\                               
  $\sigma^2$      & $0.02$ \ $(0.01)$ & $0.02$ \ $(0.01)$ & $0.03$ \ $(0.01)$ \\                               
\noalign{\smallskip}\hline                       
\end{tabular}					  		           			  
\end{center}
\end{table}

\subsection{Fluorescent lamp data}

Rosillo and Chivelet (2009) present a dataset referred here as the fluorescent lamp data. The authors analyze the lifetime of fluorescent lamps in photovoltaic systems using an analytical model whose goal is to assist in improving ballast design and extending the lifetime of fluorescent lamps. Following Rosillo and Chivelet (2009) we shall consider the nonlinear regression model (\ref{general-mixed}) with $q_i=1$, $Z_i = 0$, $\Sigma_i = \sigma^2$, ${\theta} = \left({\alpha}^{\top}, \sigma^2\right)^{\top} = \left(\alpha_0, \alpha_1, \alpha_2, \alpha_3, \sigma^2\right)^{\top}$ and
\begin{equation*}\label{nonlinear-model2}
\mu_i({\alpha}) = \frac{1}{1 + \alpha_0 + \alpha_1 x_{i1} + \alpha_2 x_{i2} + \alpha_3 x_{i2}^2}, \quad i = 1,\ldots, 14,
\end{equation*}
where the response variable is the observed lifetime/advertised lifetime ($Y$), the covariates correspond to a measure of gas discharge ($x_1$) and the observed voltage/ad- vertised voltage (measure of performance of lamp and ballast - $x_2$) and the errors are assumed to follow a normal distribution. Here we also assume a Student $t$ distribution with $\nu = 4$ degrees of freedom for the errors.

In Table 8 we present the maximum likelihood estimates, the bias-corrected estimates, the bias-reduced estimates, and the corresponding estimated standard errors. As in the previous application, the estimates present smaller standard errors under the Student $t$ model than under the normal model.

The original MLEs for $\alpha_0$ and $\alpha_3$ are bigger than the corresponding corrected and reduced versions by approximately one unit (normal and Student $t$ models). The largest differences are among the estimates of $\alpha_2$; for example, for the normal model we have $-56.33$ (MLE), $-54.45$ (bias-corrected MLE) and $-53.86$ (bias-reduced MLE). 

We now use the Akaike Information Criterion ($AIC$, Akaike, 1974), the Schwarz Bayesian criterion ($BIC$, Schwarz, 1978) and the finite sample $AIC$ ($AIC_C$, Hurvich and Tsai, 1989) to evaluate the quality of the normal and Student $t$ fits. For the normal model we have $AIC = -9.98$, $BIC = -6.79$ and $AIC_C = -2.48$. For the $t$ model we have $AIC = -11.24$, $BIC = -8.04$ and $AIC_C = -3.74$. Therefore, the $t$ model presents the best fit for this dataset, since the values of the $AIC$, $BIC$ and $AIC_C$ are smaller.

Let $$\widehat{D} = \sum_{j=1}^{14} (\widehat{Y} - \widehat{Y}_{(j)})^{\top} (\widehat{Y} - \widehat{Y}_{(j)}),$$ where $\widehat{Y}$ and $\widehat{Y}_{(j)}$ are the vectors of predicted values computed from the model fit for the whole sample and the sample without the $j$th observation, respectively. The quantity $\widehat{D}$ measures the total effect of deleting one observation in the predicted values. For a fixed sample size, it tends to be high if a single observation can highly influence the prediction of new observations. We have $\widehat{D} = 0.119, 0.120$, and $0.123$ (normal model) and $\widehat{D} = 0.101, 0.100$, and $0.095$ (Student $t$ model) when using the MLE, the bias-corrected estimate, and the bias-reduced estimate, respectively. Notice that $\widehat{D}$ is smaller for the Student $t$ model regardless of the estimate used. This is evidence that the Student $t$ model is more suitable than the normal model for predicting lifetime of fluorescent lamps in this study.

\begin{table}[!htp]
{\caption{Estimates and standard errors (given in parentheses); fluorescent lamp data.}} \label{tab13}
{
\begin{center}
\begin{tabular}{cccc} 
\hline \noalign{\smallskip}
\multicolumn{4}{c}{\mbox{Normal distribution}}  \\ 
\noalign{\smallskip}\hline\noalign{\smallskip}                                                                                      
${\theta}$  &                MLE               &          Bias-corrected MLE         &         Bias-reduced MLE            \\ 
\noalign{\smallskip}\hline\noalign{\smallskip}
$\alpha_0$ &\!\!\!\!$29.49$ \ \ \ \ \ \ \ \ \ \ \ \ \               $(5.21)$ &\!\!\!\!$28.54$ \ \ \ \ \ \ \ \ \ \ \ \ \               $( 5.66)$ &\!\!\!\!$28.25$ \ \ \ \ \ \ \ \ \ \ \ \ \               $( 5.84)$ \\                          
$\alpha_1$      & \!\!\! $9.99$ \ \ \ \ \ \ \ \ \ \ \ \ \                $(4.69)$ & \!\!\! $9.68$ \ \ \ \ \ \ \ \ \ \ \ \ \                $( 5.21)$ & \!\!\! $9.62$ \ \ \ \ \ \ \ \ \ \ \ \ \                $( 5.42)$ \\                          
$\alpha_2$      & \!\!\!\!\!\!\!\!$-56.33$ \ \ \ \ \ \ \ \ \ \ \ \ \ \ \!$(10.10)$ & \!\!\!\!\!\!\!\!$-54.45$ \ \ \ \ \ \ \ \ \ \ \ \ \ \ \!$(10.93)$ & \!\!\!\!\!\!\!\!$-53.86$ \ \ \ \ \ \ \ \ \ \ \ \ \ \ \!$(11.26)$ \\                          
$\alpha_3$      &\!\!\!\!$26.53$ \ \ \ \ \ \ \ \ \ \ \ \ \               $(4.89)$ &\!\!\!\!$25.61$ \ \ \ \ \ \ \ \ \ \ \ \ \               $( 5.28)$ &\!\!\!\!$25.31$ \ \ \ \ \ \ \ \ \ \ \ \ \               $( 5.43)$ \\                          
$\sigma^2$      &    $1.40 \times 10^{-2}$ \ \ \        $(5.00 \times 10^{-3})$ &    $1.80 \times 10^{-2}$ \ \ \        $(7.00 \times 10^{-3})$ &    $1.90 \times 10^{-2}$ \ \         $(7.00 \times 10^{-3})$ \\
  \noalign{\smallskip}\hline\noalign{\smallskip}       
\multicolumn{4}{c}{\mbox{Student $t$ distribution}}  \\ 
\noalign{\smallskip}\hline\noalign{\smallskip}                                                                                      
${\theta}$  &                MLE               &          Bias-corrected MLE         &         Bias-reduced MLE            \\ 
\noalign{\smallskip}\hline\noalign{\smallskip}
$\alpha_0$ &\!\!\!\!$30.66$ \ \ \ \ \ \ \ \ \ \ \ \ \               $(4.64)$ &\!\!\!\!$29.94$ \ \ \ \ \ \ \ \ \ \ \ \ \               $(5.05)$ & \!\!\!\!$29.85$ \ \ \ \ \ \ \ \ \ \ \ \ \               $( 5.20)$ \\         
$\alpha_1$ & \!\!\! $8.48$ \ \ \ \ \ \ \ \ \ \ \ \ \                $(4.00)$ & \!\!\! $8.24$ \ \ \ \ \ \ \ \ \ \ \ \ \                $(4.42)$ & \!\!\! $8.46$ \ \ \ \ \ \ \ \ \ \ \ \ \                $( 4.57)$ \\    
$\alpha_2$ & \!\!\!\!\!\!\!\!$-58.20$ \ \ \ \ \ \ \ \ \ \ \ \ \ \ \!$(8.94)$ & \!\!\!\!\!\!\!\!$-56.79$ \ \ \ \ \ \ \ \ \ \ \ \ \ \ \!$(9.71)$ & \!\!\!\!\!\!\!\!$-56.67$ \ \ \ \ \ \ \ \ \ \ \ \ \!$(10.00)$ \\      
$\alpha_3$ &\!\!\!\!$27.27$ \ \ \ \ \ \ \ \ \ \ \ \ \               $(4.30)$ &\!\!\!\!$26.58$ \ \ \ \ \ \ \ \ \ \ \ \ \               $(4.66)$ & \!\!\!\!$26.55$ \ \ \ \ \ \ \ \ \ \ \ \ \               $( 4.80)$ \\    
$\sigma^2$ &    $7.30 \times 10^{-3}$ \ \         $(3.60 \times 10^{-3})$ &    $9.20 \times 10^{-3}$ \ \         $(4.60 \times 10^{-3})$ &    $9.80 \times 10^{-3}$ \ \         $(4.90 \times 10^{-3})$ \\                         
\noalign{\smallskip}\hline                     
\end{tabular}	
\end{center}				  
}
\end{table}

\subsection{WHO MONICA data}

We now turn to a dataset from the WHO MONICA Project that was considered in Kulathinal et al. (2002). This dataset was first analyzed under normal distributions for the marginals of the random errors (Kulathinal et al. 2002; Patriota et al. 2009a). Thereafter, it was studied under a scale mixture of normal distributions for the marginals of the random errors (Cao et al., 2012). The approach used in the present paper is different from the others because here we consider a joint elliptical distribution for the vector of random errors. The other authors assumed that the distributions of the errors were independent, while we assume that they are uncorrelated but not independent. For our proposal, the errors will only be independent under normality.

The dataset considered here corresponds to the data collected for men ($n=38$). As describe in Kulathinal et al. (2002), the data are trends of the annual change in the event rate $(y)$ and trends of the risk scores ($x$). The risk score is defined as a linear combination of smoking status, systolic blood pressure, body mass index, and total cholesterol level. A follow-up study using proportional hazards models was employed to derive its coefficients, and provides the observed risk score and its estimated variance. Therefore, the observed response variable, $X_1$, is the average annual change in event rate (\%) and the observed covariate, $X_2$, is the observed risk score (\%). We use the heteroscedastic model (\ref{error-model-O}) with $v=m=1$ and zero covariance between the errors $\delta_{x_{1i}}$ and $\delta_{x_{2i}}$. 

Table 9 gives the MLE and the bias-corrected/reduced estimates (standard errors are given in parentheses). We considered the full sample ($n=38$) and randomly chosen sub-samples of $n = 10, 20$ and $30$ observations. 

The original MLEs for $\beta_0$, $\beta_1$ and $\mu_{x_2}$ are practically the same as their bias-corrected and bias-reduced versions for all sample sizes. The largest differences are among the estimates of ${\Sigma}_q$; for example, for $n = 10$ we have $6.17$ (MLE), $8.14$ (bias-corrected MLE) and $8.81$ (bias-reduced MLE). In general, as expected, larger sample sizes correspond to smaller standard errors. 

\begin{table}[!htp]
{\caption{Estimates and standard errors (given in parentheses); WHO MONICA data.}} 
\label{tab14} 
\begin{center}
\begin{tabular}{ccccc}
\hline\noalign{\smallskip}                                                                                          
 $n$ &  ${\theta}$ &    MLE   & Bias-corrected MLE   & Bias-reduced MLE    \\ 
\noalign{\smallskip}\hline\noalign{\smallskip}
     &    $\beta_0$      &    $-2.58$ \ $(1.34)$ &    $-2.58$ \ $(1.44)$ &    $-2.45$ \ $(1.47)$ \\                               
     &    $\beta_1$      & \ \ $0.05$ \ $(0.60)$ & \ \ $0.05$ \ $(0.63)$ & \ \ $0.07$ \ $(0.64)$ \\                               
$10$ &    $\mu_{x_2}$    &    $-1.54$ \ $(0.58)$ &    $-1.54$ \ $(0.61)$ &    $-1.53$ \ $(0.62)$ \\                               
     &    $\Sigma_{x_2}$ & \ \ $2.89$ \ $(1.50)$ & \ \ $3.22$ \ $(1.65)$ & \ \ $3.29$ \ $(1.69)$ \\                               
     &    $\Sigma_q$     & \ \ $6.17$ \ $(3.99)$ & \ \ $8.14$ \ $(4.93)$ & \ \ $8.81$ \ $(5.25)$ \\                               
\noalign{\smallskip}\hline\noalign{\smallskip}
     &    $\beta_0$      &    $-2.68$ \ $(0.65)$ &    $-2.69$ \ $(0.68)$ &    $-2.69$ \ $(0.69)$ \\                               
     &    $\beta_1$      & \ \ $0.48$ \ $(0.30)$ & \ \ $0.47$ \ $(0.31)$ & \ \ $0.43$ \ $(0.31)$ \\                               
$20$ &    $\mu_{x_2}$    &    $-1.29$ \ $(0.44)$ &    $-1.29$ \ $(0.46)$ &    $-1.29$ \ $(0.46)$ \\                               
     &    $\Sigma_{x_2}$ & \ \ $3.53$ \ $(1.25)$ & \ \ $3.73$ \ $(1.31)$ & \ \ $3.76$ \ $(1.32)$ \\                               
     &    $\Sigma_q$     & \ \ $3.00$ \ $(1.66)$ & \ \ $3.59$ \ $(1.87)$ & \ \ $3.73$ \ $(1.92)$ \\                               
\noalign{\smallskip}\hline\noalign{\smallskip}						        	 		        	 		       
     &    $\beta_0$      &    $-2.22$ \ $(0.54)$ &    $-2.22$ \ $(0.55)$ &    $-2.20$ \ $(0.55)$ \\                               
     &    $\beta_1$      & \ \ $0.43$ \ $(0.24)$ & \ \ $0.43$ \ $(0.25)$ & \ \ $0.42$ \ $(0.25)$ \\                               
$30$ &    $\mu_{x_2}$    &    $-0.77$ \ $(0.42)$ &    $-0.77$ \ $(0.42)$ &    $-0.77$ \ $(0.42)$ \\                               
     &    $\Sigma_{x_2}$ & \ \ $4.71$ \ $(1.34)$ & \ \ $4.88$ \ $(1.39)$ & \ \ $4.89$ \ $(1.39)$ \\                               
     &    $\Sigma_q$     & \ \ $4.36$ \ $(1.86)$ & \ \ $4.89$ \ $(2.01)$ & \ \ $4.88$ \ $(2.01)$ \\                               
\noalign{\smallskip}\hline\noalign{\smallskip}						        	 		        	 		       
     &    $\beta_0$      &    $-2.08$ \ $(0.53)$ &    $-2.08$ \ $(0.54)$ &    $-2.08$ \ $(0.54)$ \\                               
     &    $\beta_1$      & \ \ $0.47$ \ $(0.23)$ & \ \ $0.47$ \ $(0.24)$ & \ \ $0.46$ \ $(0.24)$ \\                               
$38$ &    $\mu_{x_2}$    &    $-1.09$ \ $(0.36)$ &    $-1.09$ \ $(0.36)$ &    $-1.09$ \ $(0.36)$ \\                               
     &    $\Sigma_{x_2}$ & \ \ $4.32$ \ $(1.10)$ & \ \ $4.44$ \ $(1.13)$ & \ \ $4.45$ \ $(1.13)$ \\                               
     &    $\Sigma_q$     & \ \ $4.89$ \ $(1.78)$ & \ \ $5.34$ \ $(1.89)$ & \ \ $5.30$ \ $(1.88)$ \\                               
\noalign{\smallskip}\hline                      
\end{tabular}		
\end{center}			  
\end{table}

\section{Concluding remarks}\label{conclusion}

\hspace{0.4cm} We studied bias correction and bias reduction for a multivariate elliptical model with a general parameterization that unifies several important models (e.g., linear and nonlinear regressions models, linear and nonlinear mixed models, errors-in-variables models, among many others). We extend the work of Patriota and Lemonte (2009) to the elliptical class of distributions defined in Lemonte and Patriota (2011). We express the second order bias vector of the maximum likelihood estimates as an weighted least-squares regression. 

As can be seen in our simulation results, corrected-bias estimators and reduced-bias estimators form a basis of asymptotic inferential procedures that have better performance than the corresponding procedures based on the original estimator. We further note that, in general, the bias-reduced estimates are less biased than the bias-corrected estimates. Computer packages that perform simple operations on matrices and vectors can be used to compute bias-corrected and bias-reduced estimates.


\section*{Appendix}\label{appendix}

{\bf Lemma A.1.} Let $z_i \sim El_{q_i}(0, \Sigma_i, g)$, and $c_{i}$ and $\psi_{i(2,1)}$ as previously defined. 
Then, 
\begin{equation*}
\begin{split}
&E\bigl(v_i{z}_i\bigr) = 0, \\
&E\bigl(v_i^2{z}_i{z}_i^{\top}\bigr) = \dfrac{4\psi_{i(2,1)}}{q_{i}}{\Sigma}_i, \\
&E\bigl(v_i^2 \vvec({z}_i{z}_i^{\top}){z}_i^{\top}\bigr) = 0, \\
&E\bigl(v_i^2\vvec({z}_i{z}_i^{\top})\vvec({z}_i{z}_i^{\top})^{\top}\bigr) = c_i\bigl(\vvec({\Sigma}_i) \vvec({\Sigma}_i)^{\top} + 2{\Sigma}_i\otimes {\Sigma}_i\bigr), \\
&E\bigl(v_i^3\vvec({z}_i{z}_i^{\top})\vvec({z}_i{z}_i^{\top})^{\top}\bigr) = - c_i^*\bigl(\vvec({\Sigma}_i) \vvec({\Sigma}_i)^{\top} + 2{\Sigma}_i\otimes {\Sigma}_i\bigr), \\
&E\bigl(v_i^3 {z}_i^{\top}A_{i(t)}{z}_i {z}_i^{\top}A_{i(s)}{z}_i{z}_i^{\top}A_{i(r)}{z}_i\bigr) = -8\widetilde{\omega}_i \bigl(\tr\{A_{i(t)}\Sigma_i\}\tr\{A_{i(s)}\Sigma_i\}\tr\{A_{i(r)}\Sigma_i\} \\ 
& \hspace{5.6cm} + 2\tr\{A_{i(t)}\Sigma_i\}\tr\{A_{i(s)}\Sigma_iA_{i(r)}\Sigma_i\} \\
& \hspace{5.6cm} + 2\tr\{A_{i(s)}\Sigma_i\}\tr\{A_{i(t)}\Sigma_iA_{i(r)}\Sigma_i\}\\
& \hspace{5.6cm} + 2\tr\{A_{i(r)}\Sigma_i\}\tr\{A_{i(t)}\Sigma_iA_{i(s)}\Sigma_i\}) \\
& \hspace{5.6cm} + 8 \tr\{A_{i(t)}\Sigma_i A_{i(s)}\Sigma_i A_{i(r)}\Sigma_i\}\bigr),
\end{split}
\end{equation*}
where $c_i^* = 8\psi_{i(3,2)}/\{q_i(q_i+2)\}$, $\psi_{i(3,2)} = E(W_g^3({r}_i){r}_i^2)$, 
$\widetilde{\omega}_i =\psi_{i(3,3)}/\{q_i(q_i+2)(q_i + 4)\}$ and $\psi_{i(3,3)} = E(W_g^3({r}_i){r}_i^3)$.
\vspace{0.2cm}

\vspace{0.5cm}
\noindent {\bf Proof: }The proof can be obtained by adapting the results of Mitchell (1989) for a matrix version.
\vspace{0.5cm}


From Lemma A.1, we can find the cumulants of the log-likelihood derivatives required to compute the second-order biases.

\vspace{0.3cm}
 
\noindent{\bf Proof of Theorem 3.1:} Following Cordeiro and Klein (1994), we write (\ref{biascorrection}) in matrix notation to obtain the second-order bias vector of $\widehat{{\theta}}$ in the form 
\begin{equation}\label{BTheta}
B_{\widehat{\theta}}(\theta) = K(\theta)^{-1}{W}\vvec(K(\theta)^{-1}),
\end{equation}
\noindent where ${W} = ({W}^{(1)},\ldots,{W}^{(p)})$ is a $p\times p^{2}$ partitioned matrix, each ${W}^{(r)}$, referring to the $r$th component of ${\theta}$, being a $p\times p$ matrix with typical $(t,s)$th element given by

\[w_{ts}^{(r)} = \frac{1}{2}\kappa_{tsr}+\kappa_{ts,r} = \kappa_{ts}^{(r)} - \frac{1}{2}\kappa_{tsr} = \frac{3}{4}\kappa_{ts}^{(r)} - \frac{1}{4}(\kappa_{t,s,r} + \kappa_{sr}^{(t)}+\kappa_{rt}^{(s)}).\]

\noindent Because $K(\theta)$ is symmetric and the $t$th element of ${W}\vvec(K(\theta)^{-1})$ is $w_{t1}^{(1)}\kappa^{1,1} + (w_{t2}^{(1)} + w_{t1}^{(2)})\kappa^{1,2} + \cdots + (w_{tr}^{(s)} + w_{ts}^{(r)})\kappa^{s,r} + \cdots   +(w_{tp}^{(p-1)} + w_{t(p-1)}^{(p)})\kappa^{p-1,p} +w_{tp}^{(p)}\kappa^{p,p}$, we may write

\begin{equation}\label{WTS}
w_{ts}^{(r)} = \frac{1}{2}(w_{tr}^{(s)} + w_{ts}^{(r)}) = \frac{1}{4}(\kappa_{ts}^{(r)}  + \kappa_{tr}^{(s)}  - \kappa_{sr}^{(t)} - \kappa_{t,s,r}).
\end{equation} 

\noindent Comparing (\ref{BTheta}) and (\ref{BIAS-vector}) we note that for the proof of this theorem it suffices to show that $F^{\top}\widetilde{{H}}\xi = W \vvec(({F}^{\top} \widetilde{H}F)^{-1})$, i.e., 
\begin{equation*}
W = F^{\top} H M H (\Phi_1,\ldots,\Phi_p).
\end{equation*}

Notice that
\begin{equation}\label{kappasr}
\begin{split}
\kappa_{sr} & = \sum_{i=1}^n \biggl\{ \frac{c_i}{2} \tr\{A_{i(r)}C_{i(s)}\} - \frac{4\psi_{i(2,1)}}{q_i}a_{i(s)}^{\top}\Sigma_i^{-1}a_{i(r)} \\
            & \hspace{1cm} -\frac{(c_i -1)}{4} \tr\{A_{i(s)}\Sigma_i\}\tr\{A_{i(r)}\Sigma_i\}  \biggr\}.
\end{split}
\end{equation}

\noindent The quantities $\psi_{i(2,1)}$ and $\psi_{i(2,2)}$ do not depend on $\theta$ and hence, the derivative of (\ref{kappasr}) with respect to $\theta_t$ is
\begin{align*}
\kappa_{sr}^{(t)} &= \sum_{i=1}^n \biggl\{ \frac{c_i}{2} \tr\{ A_{i(t)}\Sigma_iA_{i(s)}C_{i(r)} +  A_{i(s)}\Sigma_iA_{i(t)}C_{i(r)} + C_{i(ts)}A_{i(r)} \\ & + C_{i(tr)}A_{i(s)}\}\biggr\} -\sum_{i=1}^n \biggl\{ \frac{4\psi_{i(2,1)}}{q_i}\bigl(a_{i(ts)}^{\top}\Sigma_i^{-1}a_{i(r)} + a_{i(s)}^{\top}A_{i(t)}a_{i(r)} \\ 
& + a_{i(s)}^{\top}\Sigma_i^{-1}a_{i(tr)}\bigr) \biggr\}+  \sum_{i=1}^n \biggl\{ \frac{(c_i -1)}{4} \tr\{A_{i(t)}C_{i(s)} + \Sigma_i^{-1}C_{i(ts)}\}\tr\{A_{i(r)}\Sigma_i\} \biggr\}\\
& +\sum_{i=1}^n \biggl\{ \frac{(c_i -1)}{4} \tr\{A_{i(t)}C_{i(r)} + \Sigma_i^{-1}C_{i(tr)}\}\tr\{A_{i(s)}\Sigma_i\}\biggr\}.\\
\end{align*} 

\noindent Therefore,

\begin{align*}
\kappa_{st}^{(r)} + \kappa_{tr}^{(s)} - \kappa_{sr}^{(t)} &= \sum_{i=1}^n \biggl\{ \frac{c_i}{2} \tr\{ A_{i(r)}\Sigma_iA_{i(s)}C_{i(t)} +  A_{i(s)}\Sigma_iA_{i(r)}C_{i(t)} \\ 
& + 2C_{i(rs)}A_{i(t)} \}\biggr\} - \sum_{i=1}^n \biggl\{ \frac{4\psi_{i(2,1)}}{q_i}\bigl(2a_{i(t)}^{\top} \Sigma_i^{-1} a_{i(sr)} \\
& + a_{i(t)}^{\top}A_{i(s)}a_{i(r)} + a_{i(s)}^{\top}A_{i(r)}a_{i(t)} - a_{i(s)}^{\top}A_{i(t)}a_{i(r)} \bigr) \biggr\} \\
& + \sum_{i=1}^n \biggl\{ \frac{(c_i -1)}{2} \tr\{A_{i(r)}C_{i(s)} + \Sigma_i^{-1}C_{i(rs)}\}\tr\{A_{i(t)}\Sigma_i\} \biggr\}.
\end{align*} 
 
Now, the only quantity that remains to obtain is $\kappa_{t,s,r} = E(U_tU_sU_r)$. Noting that $z_i$ is independent of $z_j$ for $i\neq j$, we have 
\begin{align*}
\kappa_{t,s,r} = \frac{1}{8}\sum_{i=1}^n E\biggl\{&\bigl[ \tr\{A_{i(t)}(\Sigma_i - v_i z_iz_i^{\top}) \}\tr\{A_{i(s)}(\Sigma_i - v_i z_iz_i^{\top}) \} \\
& \tr\{A_{i(r)}(\Sigma_i - v_i z_iz_i^{\top}) \}\bigr] + 4\tr\{A_{i(t)}(\Sigma_i - v_i z_iz_i^{\top}) \}(v_i^2 a_{i(r)}^{\top}\Sigma_i^{-1} \\
& z_i z_i^{\top} \Sigma^{-1}a_{i(s)}) + 4\tr\{A_{i(r)}(\Sigma_i - v_i z_iz_i^{\top}) \}(v_i^2 a_{i(t)}^{\top}\Sigma_i^{-1}z_i z_i^{\top} \Sigma^{-1} \\
& a_{i(s)}) + 4\tr\{A_{i(s)}(\Sigma_i - v_i z_iz_i^{\top}) \}(v_i^2 a_{i(t)}^{\top}\Sigma_i^{-1}z_i z_i^{\top} \Sigma^{-1}a_{i(r)})\biggr\}.
\end{align*}
Then, by using Lemma A.1 and from (\ref{WTS}), we have, after lengthy algebra, that
\begin{equation}\label{Wr}
W^{(r)} = \sum_{i=1}^n  F_i^{\top}H_iM_{i}H_i\Phi_{i(r)},
\end{equation}
where 

\begin{equation*}
\Phi_{i(r)} = - \frac{1}{2} \left(H_i^{-1}M_i^{-1}B_{i(r)}H_iF_i + \frac{\partial F_i}{\partial \theta_r}\right),
\end{equation*}

\noindent and
\begin{align*}
B_{i(r)} &= -\frac{1}{2}
\begin{pmatrix}
\eta_{1i}  C_{i(r)}                  & &2\eta_{1i} \Sigma_i \otimes a_{i(r)}^{\top} \\
2\eta_{2i}  \Sigma_i\otimes a_{i(r)}  & & 2(c_i-1)S_{1i(r)}
\end{pmatrix} \\ \\
& -
\frac{1}{4}
\begin{pmatrix}
\eta_{1i}  \Sigma_i \tr\{C_{i(r)}\Sigma_i^{-1}\} & &2\eta_{1i} a_{i(r)}\vvec(\Sigma_i)^{\top}\\
2\eta_{1i} \vvec(\Sigma_i)a_{i(r)}^{\top}         & & 2(c_i + 8\widetilde{\omega}_i)S_{2i(r)}
\end{pmatrix},
\end{align*} with 

\begin{equation*}
\begin{split}
\eta_{1i} &= c_i^* + 4\psi_{i(2,1)}/q_i, \ \ 
\eta_{2i}  = c_i^* - 4\psi_{i(2,1)}/q_i,  \\ 
S_{1i(r)} &= \vvec(\Sigma_i)\vvec(C_{i(r)})^{\top} +\frac{1}{2}\vvec(\Sigma_i)\vvec(\Sigma_i)^{\top}\tr\{C_{i(r)}\Sigma_i^{-1}\} \ \ \mbox{and} \ \ \\
S_{2i(r)} &= \vvec(C_{i(r)})\vvec(\Sigma_i)^{\top} + \vvec(\Sigma_i)\vvec(C_{i(r)})^{\top} + 4 \Sigma_i\otimes C_{i(r)} \\
          &+ \big[\Sigma_i \otimes \Sigma_i+ \frac{1}{2}\vvec(\Sigma_i)\vvec(\Sigma_i)^{\top}\big]\tr\{C_{i(r)}\Sigma_i^{-1}\}.
\end{split}
\end{equation*}

Using (\ref{Wr}) and (\ref{BTheta}) the theorem is proved.

\vspace{0.3cm}
 
\noindent{\bf Proof of Corollary 3.1:} It follows from Theorem 3.1, eq. (\ref{BIAS-vector}), when 
\begin{equation*}
F = \mbox{block--diag}\{F_{\theta_1}, F_{\theta_2}\}, \ \widetilde{H} = \mbox{block--diag} \{\widetilde{H}_{1}, \widetilde{H}_{2}\} \ \mbox{and} \ {\xi} = (\xi_{1}^{\top}, \xi_{2}^{\top})^{\top},
\end{equation*}
\noindent where $F_{\theta_j} = \left[F_{\theta_j(1)}^{\top}, \ldots, F_{\theta_j(n)}^{\top}\right]^{\top}$ and $\widetilde{H}_{j} = \mbox{block-diag}\{\widetilde{H}_{j(1)}, \ldots, \widetilde{H}_{j(n)}\}$ for $j = 1, 2$, with $F_{\theta_1(i)} = \partial \mu_i/\partial \theta_1^{\top}$, $F_{\theta_2(i)} = \partial[\vvec(\Sigma_i)]/\partial \theta_2^{\top}$, $\widetilde{H}_{1(i)} = \frac{4\psi_{i(2,1)}}{q_i} {\Sigma}_i^{-1}$ and $\widetilde{H}_{2(i)} = c_i \left(2 \Sigma_i \otimes \Sigma_i\right)^{-1} + (c_i-1)\vvec(\Sigma_i^{-1})\vvec(\Sigma_i^{-1})^{\top}$. Furthermore, $\xi_{1} = \left[\xi_{1(1)}^\top, \ldots, \xi_{1(n)}^\top\right]^\top$ and $\xi_{2} = \left[\xi_{2(1)}^\top, \ldots, \xi_{2(n)}^\top\right]^\top$ with 
\vspace{-0.1cm}
\begin{equation*}
\begin{split}
\xi_{1(i)} =& \ - \frac{1}{2} \dot{F}_{\theta_1(i)} \ \vvec((F_{\theta_1}^{\top} \widetilde{H}_{(1)} F_{\theta_1})^{-1}), \\
\xi_{2(i)} =& \ \frac{1}{4} M_i^* P_i^* \ \vvec((F_{\theta_1}^{\top} \widetilde{H}_{(1)} F_{\theta_1})^{-1}) + \frac{1}{8} \left(M_i^* Q_i^* - 4 \ \dot{F}_{\theta_2(i)}\right) \\
            & \vvec((F_{\theta_2}^{\top} \widetilde{H}_{(2)} F_{\theta_2})^{-1}). 
\end{split}
\end{equation*} 
Also, $\dot{F}_{\theta_1(i)} = [F_{\theta_1(i)}^1, \ldots, F_{\theta_1(i)}^{p_1}]$, $\dot{F}_{\theta_2(i)} = [F_{\theta_2(i)}^1, \ldots, F_{\theta_2(i)}^{p_2}]$, $Q_i^* = [Q_{i(1)}^*, \break \ldots, Q_{i(p_2)}^*]$, $P_i^* = [P_{i(1)}^*, \ldots,$ $P_{i(p_1)}^*]$, $F_{\theta_1(i)}^r = \frac{\partial F_{\theta_1(i)}}{\partial \theta_{1(r)}}$, $F_{\theta_2(i)}^s = \frac{\partial F_{\theta_2(i)}}{\partial \theta_{2(s)}}$, where $\theta_{1(r)}$ and $\theta_{2(s)}$ are the $r$th and $s$th elements of $\theta_1$ and $\theta_2$, respectively, $r = 1, \ldots, p_1$, $s = 1, \ldots, p_2$ and

\begin{equation*}\label{Eq.MQVP}
\begin{split}
M_{i}^*  =& \ \frac{1}{c_i}\left(I_{q_i^2} - \frac{\vvec(\Sigma_i) \vvec(\Sigma_i)^\top}{2c_i + \vvec(\Sigma_i)^\top \vvec(\Sigma_i)} \right), \\
Q_{i(s)}^* =& \ \left((c_i - 1) S_{1i(s)} + \frac{1}{2} (c_i + 8\widetilde{\omega}_i)S_{2i(s)}\right) \left(\Sigma_i^{-1}\otimes \Sigma_i^{-1}\right) F_{\theta_2(i)}, \\
P_{i(r)}^* =& \ \left(2 \eta_{2i}(I_{q_i} \otimes a_{i(r)}) + \eta_{1i}\vvec(\Sigma_i)a_{i(r)}^{\top}\Sigma_i^{-1}\right) F_{\theta_1(i)}. \\
\end{split}
\end{equation*}

\section*{Acknowledgement}
We gratefully acknowledge the financial the support from CNPq and FAPESP. 

\end{document}